\documentstyle{article}

\makeatletter
\def\cicy#1(#2|#3)#4{\left(\matrix{#2}\right|\!\!
		     \left|\matrix{#3}\right)^{{#4}}_{#1}}

\def\P#1{{\bf P}^{#1}}

\def\drho#1{{\partial \; \over \partial \rho_{#1}}}

\begin{document}

\newcommand{\e}{\epsilon}
\newcommand{\w}{{\bf w}}
\newcommand{\y}{{\bf y}}
\newcommand{\z}{{\bf z}}
\newcommand{\x}{{\bf x}}
\newcommand{\N}{{\bf N}}
\newcommand{\Z}{{\bf Z}}
\newcommand{\R}{{\bf R}}
\newcommand{\Q}{{\bf Q}}
\newcommand{\C}{{\bf C}}
\newcommand{\BP}{{\bf P}}
\newcommand{\cO}{{\cal O}}
\newcommand{\cH}{{\cal H}}
\newcommand{\cM}{{\cal M}}
\newcommand{\sB}{{\sf B}}
\newcommand{\cT}{{\cal T}}
\newcommand{\cI}{{\cal I}}

\newcommand{\sE}{{\sf E}}

\newcommand{\sA}{{\sf A}}
\newcommand{\ga}{{\sf a}}
\newcommand{\es}{{\sf s}}
\newcommand{\m}{{\bf m}}
\newcommand{\cS}{{\bf S}}

\newcommand{\ihra}{\stackrel{i}{\hookrightarrow}}
\newcommand\rank{\mathop{\rm rank}\nolimits}
\newcommand\im{\mathop{\rm Im}\nolimits}
\newcommand\Li{\mathop{\rm Li}\nolimits}
\newcommand\NS{\mathop{\rm NS}\nolimits}
\newcommand\Hom{\mathop{\rm Hom}\nolimits}
\newcommand\Pic{\mathop{\rm Pic}\nolimits}
\newcommand\Hilb{\mathop{\rm Hilb}\nolimits}
\newcommand{\length}{\mathop{\rm length}\nolimits}

\newcommand\lra{\longrightarrow}
\newcommand\ra{\rightarrow}
\newcommand\cJ{{\cal J}}
\newcommand\JG{J_{\Gamma}}
\newcommand{\wvskp}{\vspace{1cm}}
\newcommand{\vskp}{\vspace{5mm}}
\newcommand{\nvskp}{\vspace{1mm}}
\newcommand{\nid}{\noindent}
\newcommand{\new}{\nvskp \nid}
\newtheorem{Assumption}{Assumption}[section]
\newtheorem{Theorem}{Theorem}[section]
\newtheorem{Lemma}{Lemma}[section]
\newtheorem{Remark}{Remark}[section]
\newtheorem{Corollary}{Corollary}[section]
\newtheorem{Conjecture}{Conjecture}[section]
\newtheorem{Proposition}{Proposition}[section]
\newtheorem{Example}{Example}[section]
\newtheorem{Definition}{Definition}[section]
\newtheorem{Question}{Question}[section]
\renewcommand{\thesubsection}{\it}

\begin{center}
\Large
{\bf On Mirror Symmetry Conjecture \\
for \\Schoen's Calabi-Yau 3 folds }

\normalsize

\vskp
Shinobu Hosono\footnote{Partly supported by Grants-in Aids for Scientific Research, the Ministry of Education, Science and Culture, Japan}

\vskp

{\it Department of Mathematics, Faculty of Science, \\ Toyama University, Toyama 930, Japan\\E-mail:hosono@sci.toyama-u.ac.jp}

\vskp
Masa-Hiko Saito\footnote{Partly supported by Grants-in Aids for Scientific Research (B-09440015), the Ministry of 
Education, Science and Culture, Japan}

\vskp
{\it Department of Mathematics, Faculty of Science, \\Kobe University, Rokko, Kobe, 657, Japan\\E-mail:mhsaito@math.s.kobe-u.ac.jp}

\vskp
{Jan Stienstra}

\vskp
{\it Department of Mathematics,   
University of Utrecht, \\Postbus 80.010, 3508 TA Utrecht, 
Netherlands\\E-mail:stien@math.ruu.nl}

\end{center}
\begin{center}
{\bf Abstracts}
\end{center}

\begin{quote}
{\footnotesize  In this paper, we  verify a part of
the Mirror Symmetry Conjecture for
Schoen's Calabi-Yau 3-fold, which is a special complete intersection
in a toric variety.
 We  calculate  a part of the
prepotential of the A-model Yukawa couplings of the Calabi-Yau
3-fold directly by means of
a theta function and Dedekind's eta function.
This gives infinitely many  Gromov-Witten
invariants, and equivalently infinitely many sets of rational curves in the
Calabi-Yau 3-fold. Using  the toric mirror construction
\cite{Batyrev-Borisov,HKTY,Sti},  we  also calculate the
prepotential of the
B-model Yukawa couplings of the mirror partner.  Comparing  the
expansion of the B-model prepotential with that of the A-model
prepotential, we check a part of the Mirror Symmetry Conjecture up to a
high order.}
\end{quote}

\vskp

\section{Introduction}
\label{intro}

Let $W$ be a generic complete intersection variety
in $\P1 \times \P2 \times \P2$ which is defined by two equations
of multi-degrees $(1, 3, 0)$ and $(1, 0, 3)$ respectively.
A generic $W$ is  a non-singular Calabi-Yau 3-fold,
which we call {\em Schoen's
Calabi-Yau 3-fold}~\cite{Sch}.
The purpose of this paper is to verify
a part of the Mirror Symmetry Conjecture
for Schoen's Calabi-Yau 3-folds.

 In \cite{COGP}, Candelas
et al. calculated the prepotential of the
B-model Yukawa couplings of the mirror of generic quintic
hypersurfaces $X$ in $\P4$ and under the mirror hypothesis they  gave
 predictions of numbers of rational curves  of
degree $d$ in $X$.  Their predictions have been  verified
mathematically only
for $d \leq 3$, that is, for numbers of lines, conics and cubic
curves (cf. e.g. \cite{E-S}).
On the B-model side one can compute as many coefficients as one wants
and thus conjecturally count curves of any degree. However it is very hard
to calculate the Gromov--Witten invariants directly on the A-model side.
 In this paper we calculate, (directly on the A-model side), a
part of the prepotential of Schoen's Calabi-Yau 3-fold $W$  which gives
infinitely many  Gromov--Witten invariants of $W$.

The main strategy of our verification is as follows:

\begin{itemize}

\item We will calculate a part of the prepotential of the A-model
Yukawa couplings (for genus zero)
of Schoen's Calabi-Yau 3-fold by using a structure of
fibration $h:W \lra \P1$ by abelian surfaces.  The theory of
Mordell-Weil lattices~\cite{Man-1,Shi-1,Saito}
allows us to calculate that part of the prepotential coming from
sections of $h$.
 Under very plausible assumptions,  we can count the ``numbers of
{\em pseudo-sections}'', which
   makes it possible for us to obtain a
very explicit description of the 1-sectional part of the A-model
prepotential (cf. Theorem~\ref{t:1-s-prep}) in 19 variables by using
a lattice
theta function for $E_8$ and Dedekind's eta function.

\item According to Batyrev-Borisov~\cite{Batyrev-Borisov} we can
construct a
mirror partner $W^*$ of $W$.
The prepotential of the B-model Yukawa couplings of $W^*$
can be defined by means of period integrals of $W^*$.
Following the recipe in \cite{HKTY,Sti} we expand the B-model
prepotential
in 3 variables by using the toric data.
These 3-variables correspond to 3 elements  of the Picard group of
$W$
coming from line bundles on the ambient space $\P1 \times \P2
\times \P2$.

\item By identifying the 3 variables with the corresponding 3
variables on the A-model side we
have two expansions which should be compared.  By a simple computer
calculation
we can verify the conjecture up to a  high order.

\end{itemize}

To state the results for the A-model side let $f_i:S_i \lra \P1$ $(i
=1,
2)$ be two generic rational elliptic surfaces.
Then Schoen's generic Calabi-Yau 3-fold
can be obtained as the fiber product $h:W =S_1 \times_{\P1}S_2 \lra
\P1$.
A general  fiber of $h$ is a product of two elliptic curves. Hence
after
fixing the zero section the
set of sections of $h$ becomes an abelian group.
In this case the group of sections $MW(W)$
is a finitely generated abelian group and admits a N\'{e}ron-Tate
height
pairing.
Let $B$ be the symmetric bilinear form associated to   the
N\'{e}ron-Tate
height pairing.
According to Shioda, we call the pair  $(MW(W), B)$ of the group and
the
symmetric bilinear form a  {\em Mordell-Weil lattice}.
Under the genericity condition for $W$ and
a suitable choice of  a N\'{e}ron-Tate height
we can easily see that  the Mordell-Weil lattice is isometric to
$E_8 \times
E_8$.  (Cf. \cite{Saito}).  There is a very
suitable set of 19 generators $[F], [L_i],
[M_j]$ $(0 \leq i, j \leq 8)$ for the Picard group of $W$.
We introduce the parameters $p, q_i, s_j$ corresponding to these
generators. The divisor class $[F]$, which  is the class of the
fiber, has
a special meaning in our context.  A homology class $\eta$ is called
$k$-sectional if the intersection number $([F], [\eta]) = k$.
Let $\Psi_A$ denote the prepotential of the A-model Yukawa couplings
of
$W$ and
$\Psi_{A, k}$ its $k$-sectional part for $k \geq 0$.
Then we have an expansion like
$$
\Psi_{A} = \mbox{topological part} + \sum_{k=0}^{\infty} \Psi_{A,
k}.
$$
Our main theorem can be stated as follows.   For detailed notation,
see
Section~\ref{s:prep}.
\begin{Theorem} $($cf. Theorem~\ref{t:mw-prep}, \ref{t:1-s-prep}$).$
Assume that Conjecture~\ref{c:cont} in
Section~\ref{s:pseudosection} holds.  Then
for a generic Schoen's Calabi-Yau 3-fold $W$ the 1-sectional
prepotential is given by
\begin{eqnarray*}
 \lefteqn{\Psi_{A,1}(p, q_0, \cdots q_8, s_0, \cdots, s_8)}
\nonumber \\
= & \sum_{n =1}^{\infty}
\frac{1}{n^3} \cdot(p \cdot \prod_{i=1}^8 (q_i \cdot s_i))^n \cdot
A(n  \tau_1, n  z_1, \ldots, n  z_8)\cdot
A(n  \tau_2, n y_1, \ldots, n  y_8).
\end{eqnarray*}
where
$$
A(\tau, x_1, \ldots, x_8) =  \Theta_{E_8}^{root}(\tau, x_1,
\ldots,x_8)
\cdot [ \sum_{m = 0}^{\infty} p(m) \cdot \exp(2 \pi m i \tau)]^{12}.
$$
Here $p(m)$ denotes the number of partitions of $m$.
\end{Theorem}
(For the definition of the various notations
see Theorem~\ref{t:mw-prep}).

Since we can prove $\Psi_{A, 0} \equiv 0$ (cf. \cite{Saito})
and $\Psi_{A,k}, k  \geq 2$ consists of terms divisible
by $p^2$, we obtain an
asymptotic expansion of $\Psi_A$ with respect to $p$:
$$
\Psi_{A} = \mbox{topological term} + p \cdot \prod_{i=1}^8(q_i s_i)\cdot  A(\tau_1, \z) \cdot
A(\tau_2, \y) + O(p^2)
$$
where $\z = ( z_1, \cdots, z_8), \y = (y_1, \cdots, y_8)$.

In Section~\ref{s:appendixII} we
show that $\Theta_{E_8}^{root}$ has a very simple expression in
terms of the
standard Jacobi theta functions.

On the other hand on the B-model side there are only 3 parameters
involved in the calculation, because the Batyrev-Borisov
construction
can only deal with  generators of the Picard group of $W$ coming
from the ambient toric variety $\P1 \times \P2 \times \P2$.
One can easily find the corresponding parameters $p = U_0, U_1 =
\exp(2 \pi i t_1),
U_2 = \exp(2 \pi i t_2)$ on the A-model side and one can
obtain the following expansion:
$$
\Psi_{A}^{res}(p, t_1, t_2) = \mbox{topological term} +
p \cdot A( 3t_1, t_1 \gamma )
\cdot A(3 t_2, t_2 \gamma) + O(p^2), 
$$
where $\gamma=(1,1,1,1,1,1,1,-1)$.
We also obtain an expansion of the $B$-model prepotential
$$
\Psi_{B}(p, t_1, t_2) = \mbox{topological term} + p \cdot B(t_1)
\cdot B(t_2) + O(p^2).
$$
Therefore in this context the mirror symmetry conjecture can be
stated as
$$
\Psi_{A}^{res}(p, t_1, t_2) \equiv \Psi_{B}(p, t_1, t_2).
$$
{From} the above  asymptotic expansion,
we come to a  concrete mathematical conjecture:
$$
A(3t, t \gamma) \equiv B(t).
$$
At this moment we can calculate both sides up to high order of
$U =  \exp(2 \pi i t)$ by using computer programs.
One can find the expansion of $A(t)$ up to
the order of 50 at the end of Section~\ref{s:rest} and also
the expansion of $B(t)$ at the end of Section~\ref{s:B-model}.

The rough plan of this paper is as follows.
In Section~\ref{s:schoen} we recall a basic property of
 Schoen's  Calabi-Yau 3-fold and its toric description.  In
Section~\ref{s:mw}
 we recall the Mordell-Weil lattice which will be essential in the
later sections.
In Section~\ref{s:prep} we first recall the definition of
Gromov-Witten invariants
and the A-model prepotential.  
We  calculate the Mordell-Weil part of the
prepotential in terms of the lattice theta function
$\Theta_{E_8}^{root}$
and state the main theorem (Theorem~\ref{t:1-s-prep}).
Section~\ref{s:pseudosection}
is devoted to counting the pseudo-sections in $W$. We also
prove Theorem~\ref{t:1-s-prep} here.
In Section~\ref{s:rest} we restrict
the parameters of the A-model prepotential in order to compare the
expansion with
that of the B-model prepotential of the mirror.  A table for
the coefficients $\{a_n\}$ of $A(3t, t \gamma)$ is given
up to order $50$(cf. Table 2 in  Section~\ref{s:rest}).
  In Section~\ref{s:B-model} after
recalling the formulation of
the mirror symmetry conjecture we calculate the B-model prepotential
of
the mirror
of Schoen's example following the recipe of \cite{HKTY,Sti}.  We
expand
the function $B(t)$  whose coefficients $\{ b_n \}$ should coincide
with
$\{ a_n \}$ if the mirror symmetry conjecture is true.  We  check
the coincidence  up to order 50.

In Appendix I (Section~\ref{s:B-model equations}) we derive the
equation of
the mirror
according to the Batyrev-Borisov
construction~\cite{Batyrev-Borisov}.
In Appendix II (Section~\ref{s:appendixII}) we give a formula for
the theta function of the lattice $E_8$.

\vskp
Let us mention some papers which are related to our work.
In the paper~\cite{D-G-W},  Donagi, Grassi and Witten
calculate the non-perturbative superpotential in $F$-theory
compactification
to four dimensions on $\P1 \times S$, where $S$ is a rational
elliptic surface. It is interesting enough to notice that the
supperpotential
in their case is also described by the lattice theta function
for $E_8$.  It is interesting that they also mention
the contribution of
Dedekind's eta function $\eta(\tau)$ to the superpotential, though
we do not know
any direct relation between $F$-theory and the  Type II theory.
  In \cite{G-P}, G\"{o}ttsche and Pandharipande
calculated the quantum cohomology of blowing-ups of $\P2$. Their
calculation for
 the blowing-up  of $9$-points in {\em general position} on $\P2$ is
certainly related to
our calculation for the rational elliptic surfaces.  Moreover, in
\cite{Y-Z}
S.-T. Yau and Zaslow describe the counting of BPS states of Type II
on K3 surfaces.
In the paper, they treated  rational curves with nodes, which may
have some relation
to our treatment of pseudo-sections.

\vskp
\section{Schoen's Calabi-Yau 3-folds}
\label{s:schoen}

Let $f_i:S_i \lra \BP^1$ ($i =1, 2$) be two smooth rational surfaces
defined over $\C$.  In this paper we always assume that an elliptic
surface has a section.

In \cite{Sch} C. Schoen showed that  the fiber product of two
rational elliptic
surfaces
$$
\begin{array}{ccccc}
  &  W= &  S_1 \times_{\BP^1} S_2  &   &   \\
  &  &   &   &   \\
    &  \swarrow p_1   & \quad     &  p_2 \searrow&    \\
 S_1 &  &  \downarrow h &   & S_2  \\
   & f_1 \searrow  &   &   \swarrow f_2   &  \\
   &   & \BP^1 &  &
 \end{array}
$$
becomes a Calabi-Yau 3-fold after  small resolutions of possible
singularities
of the fiber product.  In what follows we consider
such Calabi-Yau 3-folds which satisfy the following genericity
assumption.
\begin{Assumption}
\label{as:generic}
\begin{enumerate}
\item  The rational elliptic surfaces $f_i:S_i
\lra \BP^1$ ($i=1, 2$) are generic in the sense that the surfaces
$S_i$ are
smooth and every singular fiber of $f_i$ is
of Kodaira type $I_1$, that is, a rational curve with one node.
Then one can
see that
each  fibration $f_i$  has
 exactly 12 singular fibers of type $I_1$.  $($cf. $\cite{Kod})$.

\item Let $\Sigma_i \subset \BP^1$ be the set  of critical values of
$f_i$.
Then we assume that
$\Sigma_1 \cap \Sigma_2 = \emptyset$.

\end{enumerate}
\end{Assumption}

Under   Assumption~\ref{as:generic} the fiber product $W =  S_1
\times_{\BP^1}
S_2 $ becomes a {\em nonsingular}
Calabi-Yau 3-fold.  The following facts are well-known. (See
\cite{Kod} or
\cite{Man-1}).

\begin{Lemma}
\label{l:el}
Let $S_1, S_2, W$ be as above.
\begin{enumerate}

\item All fibers of $h:W \lra \BP^1$ have
vanishing topological Euler numbers. Hence we have $e(W) =
2(h^{1,1}(W) - h^{2,1}(W)) = 0$.

\item A generic rational elliptic surface with section
is obtained by blowing-up the 9 base points of a cubic pencil on
$\BP^2$.
 Let $\pi_1:S_1 \lra \BP^2$ and $\pi_2:S_2 \lra \BP^2$ be the
blowing-ups and $E_i, i=1,  \cdots, 9$ and $E'_j, j=1, \cdots, 9$
the
divisor classes of the
exceptional curves of $\pi_1$ and $\pi_2$ respectively.
Set $H_i = \pi_i^*(\cO_{\BP^2}(1))$.
Then we have
\begin{eqnarray}
\Pic(S_1)  &  = & \Z H_1 \oplus \Z E_1 \oplus \cdots \oplus \Z E_9,
\label{eq:pics1}\\
\Pic(S_2)  &  = & \Z H_2 \oplus \Z E'_1 \oplus \cdots \oplus \Z
E'_9.
\label{eq:pics2}
\end{eqnarray}

\item  Let $F_1$ and $F_2$ be the divisor classes of the fibers of
$f_1$ and
$f_2$ respectively.
Then we have
\begin{eqnarray}
F_1 = 3H_1  - \sum_{i=1}^{9} E_i, \quad   F_2 = 3H_2  -
\sum_{i=1}^{9} E'_i
\label{eq:f-h}
\end{eqnarray}

\item We have the following isomorphism of groups.
\begin{eqnarray}
\Pic(W) \simeq  p_1^*(\Pic(S_1)) \oplus p_2^*(\Pic(S_2))/ \Z
[p_1^*(F_1) -
p_2^*(F_2)] \label{eq:pic-w}
\end{eqnarray}

Hence the Picard number of $W$ is $h^{11}(W) = 19$. Also $h^{21}(W)
= 19$
because $e(W) = 0$.
\end{enumerate}

\end{Lemma}
{}\hfill$\Box$

We now show that Schoen's Calabi-Yau $3$-fold $W$ can also be
realized as a
complete intersection in the toric
variety  $\P1 \times \P2 \times \P2$.
Let $z_0, z_1$, $x_0, x_1, x_2$,  $y_0, y_1, y_2$
be the homogeneous coordinates of $\P1 \times \P2 \times \P2$ and
let
$$
a_0(x_0, x_1, x_2),    a_1(x_0, x_1, x_2),   b_0(y_0, y_1, y_2),
b_1(y_0, y_1,
y_2)
 $$
be generic homogeneous cubic polynomials.
Then we can assume that the generic rational elliptic surfaces
 $S_1$ and $S_2$ in Lemma \ref{l:el} are obtained
as hypersurfaces in  $\P1 \times \P2$ as in the following way.
\begin{eqnarray*}
S_1 & =&  \{P_1 = z_1\cdot a_0(x_0, x_1, x_2) - z_0 \cdot a_1(x_0,
x_1, x_2) =
0 \}  \subset \P1 \times \P2 \\
S_2 & =&  \{P_2 = z_1 \cdot b_0(y_0, y_1, y_2) - z_0 \cdot b_1(y_0,
y_1, y_2)=
0 \}  \subset \P1 \times \P2
\end{eqnarray*}
We have the natural morphisms
$$
\begin{array}{cccclc}
 & & S_i &  &\subset \P1 \times \P2 &  \\
&f_i \swarrow  &  &\searrow\pi_i & & \\
 & \P1 &   &\P2 \ ,  & &
\end{array}
$$
where $f_1 =(p_1)_{|S_i}, \pi_i = (p_2)_{|S_i}$.
Moreover,
one can easily see that $W = S_1 \times_{\BP^1} S_2 $ can be
obtained as a
complete  intersection in the toric variety $ \P1 \times \P2 \times
\P2$ of
types $(1, 3, 0) $, $(1, 0, 3)$:
$$
W  = \left\{ \begin{array}{l|c}
[z_0:z_1] \times [ x_0: x_1: x_2]\times[ y_0: y_1: y_2] &
P_1 = 0  \\
 \in  \P1 \times \P2 \times \P2  & P_2 = 0
\end{array}\right\}
$$

\section{Mordell-Weil lattices}
\label{s:mw}

The purpose of this section is a review of results on
Mordell-Weil lattices which is needed to calculate a
part of the prepotential of the A-model Yukawa couplings of Schoen's
Calabi-Yau 3-folds.  For more complete treatments the reader may
refer to
\cite{Man-1}, \cite{Shi-1}, \cite{Saito}.

We  keep the notation and assumptions of the previous section,  that
is,
let $f_i:S_i \lra \BP^1$ be  rational elliptic surfaces and  let
$h: W = S_1 \times_{\BP^1} S_2 \lra \P1$ be  the fiber product.

Let $MW(S_i), i =1, 2$ and $MW(W)$ denote the set of sections of
$f_i$ and $h$ respectively.
Since the exceptional curves of the blowing-ups
$\pi_i : S_i \lra \P2$ are the images of sections of $f_i$,
we denote by  $e_1$ and $e'_1$  the  sections
of $f_1$ and $f_2$   respectively such that
$e_1(\P1)=E_1 $ and $e'_1(\P2)= E_1'$.
 We take $e_1$ and $e'_1$ as zero sections of  $f_1$ and $f_2$
respectively. Then
$MW(S_1)$ and $MW(S_2)$ become  finitely generated abelian groups
with the
identity elements $e_1$ and $e'_1$
respectively.
  The group $MW(S_i)$ is called the Mordell-Weil group of the
rational elliptic surface
$f_i:S_i \lra \P1$.

Take the line bundles
$$
L_0 =  E_1 +  F_1 \in \Pic(S_1),\;\;  M_0 = E'_1 + F_2 \in
\Pic(S_2).
$$
Note that these line bundles are symmetric
\footnote{A line bundle on a fibration of abelian varieties
is called symmetric if it is invariant under the pull-back by
the inverse automorphism $\z \rightarrow - \z$.} and
numerically effective and $(L_0)^2 = (M_0)^2 = 1$. Hence
$L_0$ and $M_0$ are nearly ample line bundles and  $E_1$ (resp.
$E_1'$) is
the only irreducible effective curve on $S_1$ (resp. $S_2$) with
$(L_0,
E_1)_{S_1} = 0$ (resp. $(M_0, E'_1)_{S_2}= 0$).
(Here $(C, D)_{S_i}$
denotes the intersection pairing of curves on the surface $S_i$.
Later we sometimes identify this pairing with the natural
pairing $H^2(S_i) \times H_2(S_i) \lra \Z$ via Poincar\'{e} duality.
)
Thanks to Manin \cite{Man-1} we can define N\'{e}ron-Tate heights
with respect
to  $2L_0$ and $2 M_0$, that is,
quadratic forms on $MW(S_i)$ by
\begin{eqnarray}
Q_1(\sigma_1 ) = (2 L_0, \sigma_1(\P1) )_{S_1},  \quad
Q_2(\sigma_2 ) = (2 M_0,
\sigma_2(\P1) )_{S_2}
\label{eq:height}
\end{eqnarray}
for $\sigma_1 \in MW(S_1)$ and $\sigma_2 \in MW(S_2)$.

Let $B_i$ denote the positive definite symmetric bilinear form
associated to
the quadratic form $Q_i$, i.e.
$B_i(\sigma,\sigma')=
{1\over2}\{Q_i(\sigma+\sigma')-Q_i(\sigma)-Q_i(\sigma')\}.$

According to Shioda \cite{Shi-1} we call $(MW(S_i), B_i)$ {\em the
Mordell-Weil lattice}
of $f_i:S_i \lra \P1$.
Noting that our  N\'{e}ron-Tate height coincides with Shioda's
\cite{Shi-1}
we can show the following proposition.
\begin{Proposition}
Under Assumption~\ref{as:generic} in \S~\ref{s:schoen}, we have the
following isometry of
lattices.
$$
(MW(S_i), B_i) \simeq  E_8, \quad  (i =1, 2)
$$
where $E_8$ is the unique positive-definite even unimodular lattice
of rank
$8$.
\end{Proposition}
{}\hfill$\Box$

Next we consider the Mordell-Weil group $MW(W)$ of $h:W \lra \P1$,
whose
zero section corresponds to $(e_1, e'_1)$ (cf. (\ref{eq:isom})).
{From} a property of the fiber
product we have the following isomorphism:
\begin{eqnarray}
 MW(W)  &  \stackrel{\sim}{\lra} & MW(S_1) \oplus MW(S_2)
\label{eq:isom}  \\
 \sigma &  \mapsto & (\sigma_1, \sigma_2) = ( p_1 \circ \sigma, \
p_2 \circ
 \sigma) \nonumber
 \end{eqnarray}
Since the Picard group $\Pic(W)$ can be described as in
(\ref{eq:pic-w}),
we will use the following notation for the line bundles on $W$
pulled back by $p_1$ and $p_2$:
$$
\begin{array}{lll}
\ [F] =  p_1^* (F_1)=p_2^* (F_2),  &  &   \\
 \ [H_1]  =  p_1^* (H_1), &  [L_0] = p_1^* (L_0),   &  [E_i] = p_1^*
(E_i),
 \quad  (i =1, \cdots, 9),   \\
\ [H_2]= p_2^* (H_2), &  [M_0]=p_2^* (M_0), & [E_j'] =  p_2^*
(E_j'),
\quad (j =1, \cdots, 9).
\end{array}
$$
We can easily see that  $[J_0]  :=[L_0] + [M_0] $ is a symmetric
numerically effective line bundle on $W$. It  defines a
N\'{e}ron-Tate height
 on $MW(W)$ as follows.  For  each $\sigma \in MW(W)$ we set
\begin{eqnarray}
Q_W(\sigma) := ([2J_0], [\sigma(\P1)])_W.
\end{eqnarray}
Here $( \quad , \quad  )_W$ denotes the natural pairing $H^2(W)
\times H_2(W)
\rightarrow \Z$.  Note that
the zero section of  $MW(W)$ is  $0_W =(e_1, e'_1) $ and $Q_W(0_W) =
0$.
{From} this we obtain the Mordell-Weil lattice $(MW(W),
B_W)$ where $B_W$ denotes
the symmetric bilinear form  associated to $Q_W$.
Moreover we obtain the following relation for each section $\sigma
\in MW(W)$:
\begin{eqnarray}
Q_W(\sigma)& =& ([2J_0], [\sigma(\P1)])_W  \nonumber \\
&=&  (2L_0, [\sigma_1(\P1)])_{S_1} +  (2M_0, [\sigma_2(\P1)])_{S_2}
= Q_1(\sigma_1) + Q_2(\sigma_2). \label{eq:qw}
\end{eqnarray}

Therefore we find the following

\begin{Proposition}
The N\'{e}ron-Tate height with respect to $[2J_0]$ on $MW(W)$ gives
a
 lattice structure on $MW(W)$ which induces the isometry:
$$
(MW(W), B_W)  \simeq (MW(S_1), B_1) \oplus (MW(S_2), B_2)  \simeq
E_8 \oplus
E_8.
$$

\end{Proposition}
{}\hfill$\Box$

There are  natural maps
$$
\begin{array}{cccl}
j :&MW(S_i) &\lra  &  H_2(S_i) \\
      & \sigma & \mapsto & j(\sigma) = [\sigma(\P1)]=
       \mbox{the homology class of the curve } \sigma(\P1)    \\
    &  &  &   \\
j :&MW(W) &\lra  &  H_2(W) \\
      & \sigma & \mapsto & j(\sigma) = [\sigma(\P1)].
\end{array}
$$
For each section $\sigma_i \in MW(S_i)$, we can always find a
birational
morphism $\varphi_i:S_i \lra \overline{S_i}$ which contracts only
the image of
the
section $\sigma_i(\P1)$.  This implies the  following lemma.

\begin{Lemma}\label{l:j-inj}
The maps $j$ are injective.
\end{Lemma}
{}\hfill$\Box$

Note that the maps  $j$ are  {\em not}  homomorphisms of
groups.\footnote{However, Shioda~\cite{Shi-1} obtained a way to
modify the map
$j$ to obtain a natural homomorphism.  See \cite{Shi-1} or
\cite{Saito}.}

\vskp
Next we will choose other generators of $\Pic(S_i)$.
These generators will be used for defining the parameters in
which we will expand  the prepotential of the A-model Yukawa
coupling of
Schoen's Calabi-Yau 3-folds.
   Let $(MW(S_i), B_i)$
be the Mordell-Weil lattices of $S_i$, which are isometric to the
lattice
$E_8$.
We choose a set of simple roots $ \{ \alpha_1, \alpha_2,  \cdots
\alpha_8 \} $
of
$E_8$
whose intersection pairing will be given by the following Dynkin
diagram.

\begin{picture}(315,120)(-50,0)

\put(40,70){\circle{30}}
\put(35,68){$\alpha_1$}
\put(55,70){\line(1,0){14}}
\put(85,70){\circle{30}}
\put(81,68){$\alpha_2$}
\put(100,70){\line(1,0){14}}
\put(130,70){\circle{30}}
\put(125,68){$\alpha_3$}
\put(145,70){\line(1,0){14}}
\put(130,25){\circle{30}}
\put(130,54.5){\line(0, -1){14}}
\put(125,23){$\alpha_8$}
\put(175,70){\circle{30}}
\put(170,68){$\alpha_4$}
\put(190,70){\line(1,0){14}}
\put(220,70){\circle{30}}
\put(215,68){$\alpha_5$}
\put(235,70){\line(1,0){14}}
\put(265,70){\circle{30}}
\put(260,68){$\alpha_6$}
\put(280,70){\line(1,0){14}}
\put(310,70){\circle{30}}
\put(306,68){$\alpha_7$}

\put(110,-15){Figure 1.}

\end{picture}
\vskp
\vskp

We also choose $a_1, \cdots, a_8 \in MW(S_1)$ (resp. $b_1, \cdots,
b_8 \in
MW(S_2)$)
corresponding with the roots of $MW(S_1)$ (resp. $MW(S_2)$)
so that the numbering of the
roots is the same as in Figure 1.

For each section $\sigma \in MW(S_i)$, one can define a
translation automorphism $T_{\sigma}:S_i \lra S_i$:
$$
\begin{array}{ccccc}
S_i &  &  \stackrel{T_{\sigma}}{\lra} &  &  S_i \\
 & f_i \searrow  &   &   \swarrow  f_i  &  \\
   &   & \BP^1. &  &
 \end{array} \
$$

Pulling back the line bundles $L_0$ and $M_0$ by
the translation automorphisms $T_{a_i}$ and $T_{b_j}$ respectively,
we define the line bundles
\begin{eqnarray}
L_i = T^*_{a_i} (L_0) \in \Pic(S_1), \quad M_j = T^*_{b_j} (M_0) \in
\Pic(S_2),
\end{eqnarray}
for $1 \leq i, j \leq 8$.

Then  for each section $\sigma_i \in MW(S_i)$ we have
\begin{eqnarray*}
(L_i, j(\sigma_1))_{S_1} &= &(T_{a_i}^*(L_0), j(\sigma_1))_{S_1}
=(L_0, j(\sigma_1 + a_i))_{S_1} =
\frac{1}{2}Q_1(\sigma + a_i)
\\
(M_i, j(\sigma_2))_{S_2} &= &(T_{b_i}^*(M_0), j(\sigma_2))_{S_2}
=(M_0, j(\sigma_2 + b_i))_{S_2} =
\frac{1}{2}Q_2(\sigma_2 + b_i)
\end{eqnarray*}

Now it is easy to see the following:

\begin{Lemma}\label{l:gen}
\begin{enumerate}
\item  The classes $F_1, L_0, L_1, \cdots, L_8 $  $($
resp. $F_2, M_0, M_1, \cdots, M_8   )$ are generators
of $\Pic(S_1)$ (resp. $\Pic(S_2)$).

\item  $\Pic(W) $ is generated by
$[F], [L_0], \cdots, [L_8], [M_0], [M_1], \cdots, [M_8]$.

\item  For $\sigma  \in MW(W)$ set $\sigma_i = \sigma \circ p_i$.
Then  we
have
$$
\begin{array}{lcl}
([F], j(\sigma))_W & =& 1  \\[.5em]
([L_i], j(\sigma))_W & =& \frac{1}{2} Q_1(\sigma_1 + a_i), \\[.5em]
([M_i], j(\sigma))_W & =& \frac{1}{2} Q_2(\sigma_2 + b_i).
\end{array}
$$
\end{enumerate}
\end{Lemma}
{}\hfill$\Box$

Moreover, in order to see the relation between the A-model and
the B-model later,
we have to express  $H_1$ and $H_2$ by
$F_1, \{ L_i \}$ and $F_2, \{ M_j \}$.
Obviously, we only have to see the case of $H_1$.
Recall that the exceptional curves $\{ E_i \}$ in (\ref{eq:pics1})
are the images of sections of $f_1$.
We denote by $e_i \in MW(S_1)$ the section corresponding to $E_i$;
hence we
have $e_i(\P1) = E_i$.  In particular, $e_1$ is the zero section of
$f_1:S_1 \lra \P1$.
As for the system of roots, one can take the following elements:
$$
a_1 = e_2, \ a_2 = e_3 - e_2, \ a_3 = e_4 - e_3,  \cdots, a_7 = e_8
- e_7,
$$
and
$$
a_8 =  e_2 + e_3 - \frac{1}{3} \sum_{i =2}^9 e_i.
$$
Here all sums are taken in the  Mordell-Weil group.
We denote by $(\sigma) \in H^2(S_1, \Z)$
 the divisor class of the curve $\sigma(\P1) \subset S_1$.
Since $L_0 = E_1 +  F_1 =  (e_1) + F_1$, we see that
$$
L_i = T^*_{a_i}(L_0) =  T^*_{a_i}((e_1) + F_1) = (-a_i) + F_1.
$$
Moreover we can see the following relation. (For divisor classes
$(-a_i)$ one
may
refer to \cite{Saito}).
$$
\begin{array}{lclcl}
L_0  & = & E_1 + F_1& &    \\
L_1 & = & (-a_1) + F_1 & = & 2 E_1 -E_2 + 3 F_1 \\
L_2 & = & (-a_2) + F_1 & = & E_1 + E_2 - E_3 + 2 F_1   \\
L_3 & = & (-a_3) + F_1 & = & E_1 + E_3 -E_4 + 2 F_1 \\
L_4 & = & (-a_4) + F_1 & = & E_1 + E_4 - E_5 + 2F_1 \\
L_5 & = & (-a_5) + F_1 & = &  E_1 + E_5 - E_6 + 2 F_1\\
L_6 & = & (-a_6) + F_1 & = &  E_1 + E_6 - E_7 + 2 F_1\\
L_7 & = & (-a_7) + F_1 & = &  E_1 + E_7 - E_8 + 2 F_1\\
L_8 & = & (-a_8) + F_1 & = &
 \frac{1}{3}\sum_{i=1}^9 E_i - (E_2 + E_3) +\frac{4}{3} F_1.
\end{array}
$$
Recall also the relation~(\ref{eq:f-h}):
\begin{eqnarray*}
F_1 = 3H_1  - \sum_{i=1}^{9} E_i.
\end{eqnarray*}
{From} these linear relations one easily derives the following:
\begin{Lemma}
One  has the following relation in $\Pic(S_1)$:
\begin{eqnarray}
H_1 = 2 F_1 + 5 L_0 -2 L_1 - L_2   + L_8.  \label{eq:h-f}
\end{eqnarray}
\end{Lemma}
{}\hfill$\Box$

\section{The prepotential of the A-model Yukawa couplings and its
1-sectional
part}
 \label{s:prep}

In this section
we  summarize a result in (\cite{Saito}) on the Mordell-Weil part
of the prepotential of the
A-model Yukawa coupling of Schoen's Calabi-Yau 3-folds.
The main theorems are
Theorem~\ref{t:mw-prep} and Theorem~\ref{t:1-s-prep}.

Following  Section 3.3 in \cite{Mor-Math}, we define
the (full) {\em A-model Yukawa coupling}
for a Calabi-Yau 3-fold $X$  as a triple product on  $H^2(X, \Z)$:
\begin{eqnarray}
\Phi_A(M_1, M_2, M_3) = (M_1, M_2, M_3) + \sum_{0 \not\equiv
\eta \in H_2(X, \Z)}
\Phi_{\eta}(M_1, M_2, M_3)  \frac{q^{\eta}}{1 - q^{\eta}}.
\label{eq:a-yukawa}
\end{eqnarray}
Here $M_1, M_2, M_3$  are elements of
$ H^2(X, \Z) \cong \Pic(X)$
and
$\Phi_{\eta}(M_1, M_2, M_3)$
is  the Gromov-Witten Invariant for $ \eta \in H_2(X, \Z)$ and
$M_1, M_2, M_3$.

Moreover, we have (cf. Section 3.2, \cite{Mor-Math}):
\begin{eqnarray}
\Phi_{\eta}(M_1, M_2, M_3) &= &n(\eta) (M_1, \eta)(M_2,
\eta)(M_3, \eta).  \label{eq:g-w.inv}
\end{eqnarray}

Here $(M_i,  \eta)$ denote the natural pairing of $M_i \in H^2(X)$
and $\eta \in H_2(X) $ and
$n(\eta)$ denotes the number of simple rational curves
$\varphi:\BP^1
\lra  X $ with $\varphi_*([\BP^1]) = \eta$.
A more precise definition by $J$-holomorphic curves can
be found in  \cite{McD-S-1} and Lecture 3 of \cite{Mor-Math}.

\vskp
The full Yukawa coupling $\Phi_A$ has the {\em prepotential}
$\Psi_A$ defined
by
\begin{equation}
\Psi_A =  \mbox{(topological term)} + \sum_{0 \not\equiv \eta
\in H_2(X, \Z)}  n(\eta) \Li_3( q^{\eta}),
\end{equation}
where
\begin{equation}
\Li_3( x ) = \sum_{n=1}^{\infty} \frac{x^n}{n^3}
\end{equation}
is the trilogarithm function.

In general it is very difficult to calculate the prepotential
of the A-model Yukawa coupling. Even for Schoen's Calabi-Yau 3-fold,
we can not calculate the full prepotential, but by using the
structure of its
Mordell-Weil lattice,  we can calculate a part of the
prepotential $\Psi_A$ whose summation is taken just over
the homology 2-cycles of  sections of $h:W \lra \P1$.
Later we will  extend the  summation to all
homology classes of {\em pseudo-sections }
(see Section~\ref{s:pseudosection}). (Cf. \cite{Saito}).

\begin{Definition} \label{def:prep}
{\rm For  Schoen's generic Calabi-Yau 3-fold $W$
 we define the {\em
 Mordell-Weil part of the prepotential of the A-model Yukawa
 coupling}  by
\begin{eqnarray}
\Psi_{A, MW(W)} = \sum_{\sigma \in MW(X)} n(j(\sigma))
\Li_3 ( q^{j(\sigma)}). \label{eq:prepdef}
\end{eqnarray}
Here again $j(\sigma) $ denotes the homology class of the
image $\sigma(\P1)$. }
\end{Definition}

\begin{Definition}
\label{def:sec-prep}
 {\rm   We define the
$k$-sectional part of the prepotential of the A-model Yukawa coupling by
\begin{eqnarray}
\Psi_{A, k} = \sum_{0 \not\equiv \eta
\in H_2(X, \Z), \ (F, \eta)= k}  n(\eta) \Li_3( q^{\eta}).
\label{eq:sec-prep}
\end{eqnarray}
 Recall that we denote by $[F]$
the divisor class of the fiber of $h:W \lra \P1$.}
\end{Definition}

Obviously, we have the expansion
\begin{eqnarray}
\Psi_{A}  = \mbox{topological term} + \sum_{k = 0}^{\infty} \Psi_{A,
k}.
\label{eq:expansion}
\end{eqnarray}
We are interested in calculating the functions 
$\Psi_{A, MW(W)}$ and $\Psi_{A, 1}$.

\begin{Remark}
{\rm We will find a difference in
$\Psi_{A, MW(W)}$ and $ \Psi_{A, 1}$, which will be explained
in the next section by introducing the notion {\it pseudo-section}.
}
\end{Remark}

We first recall a result in \cite{Saito} on the calculation of
$\Psi_{A, MW(W)}$ by using the theta function of the Mordell-Weil
lattice.  We need to introduce the special formal parameters in
order
to get explicit expansions of $\Psi_{A, MW(W)}$.

Let $f_i:S_i \lra \BP^1$ be two generic rational elliptic surfaces
 and let $h:W \lra \P1$ be the Calabi-Yau
3-fold as in Section~\ref{s:schoen}.  Then from Lemma~\ref{l:gen}
$\Pic(W) $ is generated by $[F], [L_0], \cdots, [L_8],$
$[M_0], [M_1], \cdots, [M_8]$.
We introduce formal parameters  $p, q_i, s_j$  for
$0 \leq i, j \leq 8 $ corresponding to these
generators:
\begin{eqnarray}
[F] \leftrightarrow p, \quad [L_i] \leftrightarrow q_i, \quad
[M_j] \leftrightarrow s_j.  \label{eq:parameter}
\end{eqnarray}

By using the formal parameters we can associate to $\sigma \in
MW(W)$ the monomials
\begin{eqnarray}
q^{\sigma} = \prod_{i=0}^8 q_i^{([L_i], j(\sigma))_W}, \quad
s^{\sigma} = \prod_{i=0}^8 s_i^{([M_i], j(\sigma))_W}.
\label{eq:prod}
\end{eqnarray}
and
\begin{equation}
T^{\sigma} = p^{([F], j(\sigma))_W} \cdot q^{\sigma} \cdot
s^{\sigma}
= p \cdot q^{\sigma} \cdot s^{\sigma}. \label{eq:tsigma}
\end{equation}
Here $( \quad, \quad)_W:H^2(W) \times H_2(W) \lra \Z $ is
the natural pairing.
Note that all line bundles $[F], [L_i], [M_j]$
are numerically effective. Hence all  exponents in $T^{\sigma}$
are non-negative.  Now we can expand $\Psi_{A, MW(W)}$ in
the parameters $p, q_i, s_j$.

 \begin{Theorem}\label{t:mw-prep}
\begin{eqnarray}
\lefteqn{\Psi_{A, MW(W)}(p, q_0, \cdots q_8, s_0, \cdots, s_8)}
\nonumber \\
= & \sum_{n =1}^{\infty}
\frac{1}{n^3} \cdot (p \cdot \prod_{i=1}^8 (q_i \cdot s_i))^n \cdot
\Theta_{E_8}^{root}(n  \tau_1, n \cdot \z) \cdot
\Theta_{E_8}^{root}(n  \tau_2, n \cdot \y)
\label{eq:main}
\end{eqnarray}
Here, we set
\begin{eqnarray*}
\z = (z_1, \ldots, z_8), \quad \y = (y_1, \ldots, y_8),
\end{eqnarray*}
\begin{eqnarray}
\exp( 2 \pi i \tau_1) = \prod_{i=0}^8 q_i, &  \exp(2 \pi i z_i)=
q_i \quad \mbox{for }  1 \leq i \leq 8       \\
\exp(2 \pi i \tau_2) = \prod_{i=0}^8 s_i,   &  \exp(2 \pi i y_j)=
s_j    \quad  \mbox{ for }  1 \leq j \leq 8
\end{eqnarray}
and
\begin{eqnarray}
\Theta_{E_8}^{root}(\tau, z_1, \cdots, z_8) =  \sum_{\gamma \in E_8}
\exp(2 \pi i ((\tau/2) Q(\gamma) + B(\gamma, \sum_{j =1}^{8} z_j
\alpha_j)),
\label{eq:roottheta}
\end{eqnarray}
where  $\{\alpha_1, \cdots, \alpha_8 \}$ is the root system of $E_8$
as in
Figure 1 and $B$ is the symmetric bilinear form on $E_8$.
\end{Theorem}

The following lemma is easy but essential to calculate the
prepotential.
\begin{Lemma} \label{l:mult=1}
For each section $\sigma \in MW(W)$ the contribution
of the homology $2$-cycle
$j(\sigma) =[\sigma(\P1)]$
to the Gromov-Witten invariant~(\ref{eq:g-w.inv})  is one, that
is, $ n(j(\sigma)) =  1 $
\end{Lemma}

\noindent
{\it Proof.}  According to Lemma~\ref{l:j-inj} $MW(W)$ can be
considered
as
a subset of $H_2(W, \Z)$  via the map $j$. Moreover the rational
curve
$C = \sigma(\P1) \subset W$ has the
normal bundle $\cO_{\P1}(-1) \oplus \cO_{\P1}(-1)$. Hence we have
$n(j(\sigma)) =1$.
 {}\hfill$\Box$

\vskp
\noindent
{\it Proof of Theorem~\ref{t:mw-prep}.} Recalling  the
isomorphism~(\ref{eq:isom}),
one can write $\sigma \in MW(W)$ as $(\sigma_1, \sigma_2) \in
MW(S_1) \oplus MW(S_2) \simeq E_8 \oplus E_8$.  Since
$Q_1(a_i) = Q_2(b_j) = 2$  for
$1 \leq i, j \leq 8$ we obtain from Lemma~\ref{l:gen}
\begin{eqnarray}
([L_i], [\sigma(\P1)])_W  =  \frac{1}{2} Q_1(\sigma_1
+ a_i) = \frac{1}{2} Q_1(\sigma_1) + B_1(\sigma_1, a_i) + 1,  \\
\quad ([M_j], [\sigma(\P1)])_W = \frac{1}{2} Q_2(\sigma_2
+ b_j) = \frac{1}{2} Q_2(\sigma_2) + B_2(\sigma_2, b_j) + 1.
\end{eqnarray}
Therefore one has
\begin{eqnarray*}
\lefteqn{ q^{\sigma} = \prod_{i=0}^{8}(q_i)^{(1/2) Q_1(\sigma_1 +
a_i) } } \\
& &  =  (\prod_{i=0}^{8}(q_i))^{(1/2) Q_1(\sigma_1)} \cdot
(\prod_{i=1}^{8}
(q_i)^{B_1(\sigma_1, a_i)}) \cdot (\prod_{i=1}^{8}q_i)  \\
& & =(\prod_{i=1}^{8}q_i) \cdot \exp(2 \pi i ((1/2) Q_1(\sigma_1)
\tau_1 +
\sum_{i=1}^8 z_i B(\sigma_1,  a_i) ),
\end{eqnarray*}
and a similar expression for $s^{\sigma}$. {From} these formulas one
can
obtain
\begin{eqnarray}
\lefteqn{ (T^{\sigma})^n = (p \cdot q^{\sigma} \cdot s^{\sigma})^n
}  \nonumber \\
\lefteqn{=( p \prod_{i=1}^8 (q_i \cdot s_i) )^n \times }
\label{eq:tsigma1}\\
  & \times \exp(2 \pi i n ((\tau_1/2) Q_1(\sigma_1)  +
B_1(\sigma_1, \z ) +
  (\tau_2/2) Q_2(\sigma_2)  +  B_2(\sigma_2, \y ) ) \nonumber
\end{eqnarray}
where we set $\z = \sum_{i=1}^8 z_i a_i$ and $\y = \sum_{i=1}^8 y_i
b_i$.
Therefore, if we take the summation of $(T^{\sigma})^n$ over $\sigma
=
(\sigma_1, \sigma_2)
\in E_8 \oplus E_8$, we obtain the following formula:
\begin{eqnarray}
\lefteqn{\sum_{\sigma \in MW(W)} (T^{\sigma})^n}  \nonumber \\
&= & ( p \prod_{i=1}^8 (q_i \cdot s_i) )^n \cdot
\Theta_{E_8}^{root}(n  \tau_1, n \cdot z ) \cdot
\Theta_{E_8}^{root}(n \tau_2, n \cdot \y)
\end{eqnarray}

Now thanks to Lemma~\ref{l:mult=1},
 we can  calculate the prepotential as follows:
$$
\begin{array}{lcl}
\Psi_{A, MW(W)}& = & \sum_{\sigma \in MW(W)} \Li_3( T^{\sigma})  \\[.5em]
	 & = & \sum_{\sigma \in MW(W)} ( \sum_{n =
	 1}^{\infty}\frac{(T^{\sigma})^n}{n^3} ) \\[.5em]
  &= & \sum_{n = 1}^{\infty} \frac{1}{n^3} \times [\sum_{\sigma \in
MW(W)}
  (T^{\sigma})^n) ] \\[.6em]
 & = &  \sum_{n = 1}^{\infty}
\frac{1}{n^3} \cdot ( p \prod_{i=1}^8 (q_i \cdot s_i))^n \times \\[.5em]
& & \hspace{1cm} \times \Theta_{E_8}^{root}(n  \tau_1, n \cdot \z)
\cdot
\Theta_{E_8}^{root}(n \tau_2, n \cdot \y).
\end{array}
$$
This completes the proof of Theorem~\ref{t:mw-prep}.
{}\hfill$\Box$

\vskp

For the 1-sectional part $\Psi_{A, 1}$
of the prepotential, we can show the following theorem, whose
proof  can be found in  Section~\ref{s:pseudosection}.

\begin{Theorem}
\label{t:1-s-prep}  Assume that Conjecture~\ref{c:cont} in
Section~\ref{s:pseudosection} holds.  Then,  under the same notation
as in Theorem~\ref{t:mw-prep},
for a generic Schoen's Calabi-Yau 3-fold $W$ the 1-sectional
prepotential is given by
\begin{eqnarray}
 \lefteqn{\Psi_{A,1}(p, q_0, \cdots q_8, s_0, \cdots, s_8)}
 \\
= & \sum_{n =1}^{\infty}
\frac{1}{n^3}\cdot (p \cdot \prod_{i=1}^8 (q_i \cdot s_i))^n \cdot
A(n  \tau_1, n \cdot \z)\cdot
A(n  \tau_2, n \cdot \y),
\nonumber
\label{eq:1-s-prep}
\end{eqnarray}
where
\begin{eqnarray}
A(\tau, \x) & = &  \Theta_{E_8}^{root}(\tau, \x )
\cdot [ \sum_{m = 0}^{\infty} p(m) \cdot \exp(2 \pi m i \tau)]^{12}.
\\
      & = & \Theta_{E_8}^{root}(\tau, \x ) \cdot
[\frac{1}{\prod_{n\geq1}( 1-
      \exp(2 \pi n i \tau))} ]^{12}
\end{eqnarray}
Here $p(m)$ denotes the number of partitions of $m$.
\end{Theorem}

\begin{Remark}\label{rem:root}
{\rm In order to identify the theta function
$\Theta_{E_8}^{root}(\tau, \z)$ in
(\ref{eq:roottheta})
with the theta function $\Theta_{E_8}(\tau, \w)$ of (\ref{eq:e8}) in
Appendix II we should apply the linear transformation $\w \lra \z$,
for
$\w = \sum_{i=1}^8 w_i \e_i$ and $\z = \sum_{i=1}^8 z_i \alpha_i$.
Fix an embedding $E_8 \subset \R^8$, that is,
$\alpha_i $ should have coordinates in $\R^8$.
For example, we can choose
$$
\alpha_1  =   \frac{1}{2}(\e_1 + \e_8) -
\frac{1}{2}(\e_2 + \e_3 + \e_4 + \e_5 + \e_6 + \e_7) \hspace{1.2cm}
$$
$$
\begin{array}{cllll}
\hspace{1cm} &\alpha_2 = \e_2 - \e_1 &  \alpha_3 = \e_3- \e_2 &
\alpha_4 = \e_4 - \e_3  & \alpha_5  =  \e_5 - \e_4\\
\hspace{1cm}  & \alpha_6  =  \e_6 - \e_5 & \alpha_7 = \e_7 - \e_6
&   \alpha_8  =  \e_1 + \e_2&
\end{array}
$$
(Note that the  numbering of roots is the same as in Figure 1.)}
\end{Remark}
 \begin{Remark}
 \label{r:first}
 {\rm
 In expansion~(\ref{eq:expansion}), we see that
 each term of the expansion of $\Psi_{A, k}$ for
 $k \geq 2$
is divisible by $p^2$. Moreover we can see that
$\Psi_{A,0} \equiv 0$. (Cf. \cite{Saito}).
 Therefore, Theorem~\ref{t:1-s-prep} shows that
if we expand the full A-model prepotential $\Psi_A$ in the variables
in
 $p, q_i, s_j$,  we have the following expansion:
 \begin{eqnarray}
 \lefteqn{\Psi_{A}(p, q_0, \cdots q_8, s_0, \cdots, s_8)} \nonumber
\\
  &= & \mbox{topological term}
 +  (p \cdot \prod_{i=1}^8 (q_i \cdot s_i)) \cdot
A(  \tau_1,  \z )\cdot
A(  \tau_2,  \y) + O(p^2).
\label{eq:1prep-exp}
\end{eqnarray}
 }
 \end{Remark}

\section{Counting Pseudo-Sections and Proof of Theorem 4.2}
 \label{s:pseudosection}

In Section~\ref{s:prep}, we see  differences between the
two prepotentials $\Psi_{A, MW(W)}$ and $\Psi_{A, 1}$.
Looking at the formulas~(\ref{eq:main}) and (\ref{eq:1-s-prep}) one
can observe
that $\Psi_{A, MW(W)}$ and $\Psi_{A, 1}$ are essentially produced by
the
functions
\begin{eqnarray}
 \Psi_{A, MW(W)}  & \leftrightarrow   &  \Theta_{E_8}^{root}(\tau,
\x)
 \label{eq:theta}\\
\Psi_{A, 1} & \leftrightarrow  & A(\tau, \x)
=  \Theta_{E_8}^{root}(\tau, \x)
 \cdot [ \sum_{m = 0}^{\infty} p(m) \cdot \exp(2 \pi m i \tau)]^{12}
\label{eq:A}
\end{eqnarray}
As we see in Section~\ref{s:prep} the geometric meaning of the
function
$ \Theta_{E_8}^{root}$ is clear, that is, it is the generating
function of the contributions of pure sections of $h:W \lra \P1$.
However, the meaning of the factor
$$
[ \sum_{m = 0}^{\infty} p(m) \cdot \exp(2 \pi m i \tau)]^{12} =
\exp(\pi i \tau)\cdot \eta(\tau)^{-12}
$$
was  mysterious at least in the geometric sense.\footnote{
The similar  factor are also discussed  in the papers \cite{D-G-W} and
\cite{Y-Z}.}
In this section, we give a geometric explanation of this
factor assuming one very plausible Conjecture~\ref{c:cont}, and we
give a
proof of Theorem~\ref{t:1-s-prep}.
Our answer seems to be very  simple and natural at
least in a mathematical sense.

For this purpose we give the following:

\begin{Definition}
{\rm We call a 1-dimensional connected subscheme  $C$  of $W$  a
{\em
pseudo-section}  if
$C \subset W$ has no embedded  point and
\begin{eqnarray}
([F], C)_W = 1,
\label{eq:p-section}
\end{eqnarray}
and the normalization  $\tilde{C}_{red}$ of the reduced structure
$C_{red}$ is a
sum of $\P1$s.}
\end{Definition}
\begin{Example}{\rm
The image  $\sigma(\P1)$ of a section $\sigma \in MW(W)$ is  a
pseudo-section.}
\end{Example}
\begin{Example}{\rm
Both rational elliptic surfaces $f_i:S_i \lra \P1 (i =1, 2)$
have 12 singular fibers of type $I_1$ (in Kodaira's
notation~\cite{Kod}):
\begin{eqnarray}
D_1, D_2, \cdots, D_{12} \subset S_1, \\
D'_1, D'_2, \cdots, D'_{12} \subset S_2.
\end{eqnarray}
We set $d_i = f_1(D_i) \in \P1$
and $ d'_i = f_2(D'_i) \in \P1$, the supports of the
singular fibers.  By Assumption~\ref{as:generic}  in
Section~\ref{s:schoen},
the points $d_1, \cdots, d_{12}, d'_1, \cdots, d'_{12} $ are
distinct on
$\P1$.
Take any  $\sigma \in MW(W)$ and set $\sigma_1 =p_1 \circ \sigma \in
MW(S_1),
\sigma_2 = p_2 \circ \sigma
\in MW(S_2)$.
 Hence $\sigma(\P1) \subset W_{|\sigma_2(\P1)}$.
Now we take a singular fiber $D_1  \subset W_{|\sigma_2(\P1)} \simeq
S_1$, then
$$
\sigma(\P1) + D_1 \subset W_{|\sigma_2(\P1)} \subset W
$$
is a pseudo-section.  Since we have
\begin{eqnarray}
([L_i], D_1)_W &=& (L_i, D_1)_{S_1} = (L_i, F_1)_{S_1} = 1
\label{eq:p-degree1} \\
([M_j], D_1)_W &=& (M_j, \sigma_2(d_1))_{S_2}  = 0,
\quad   \label{eq:p-degree2}
\end{eqnarray}
we obtain
\begin{eqnarray}
 ([L_i], \sigma(\P1)+ D_1)_W = ([L_i], \sigma(\P1))_W + 1, \\
 ([M_j],\sigma(\P1)+ D_1)_W = ([M_j], \sigma(\P1))_W.
\end{eqnarray}
 }
\end{Example}

\begin{Example} {\rm
More generally, to a pure section of $h:W \lra \P1$
we can add  many rational curves coming from singular
fibers of type  $I_1$ and
also with multiplicity.
Fix a section $\sigma \in MW(W)$ and set $\sigma_1 =p_1 \circ
\sigma$, $\sigma_2
= p_2 \circ \sigma$ as before.
Consider the following (reduced) closed points:
\begin{eqnarray}
\sigma_1(d'_i)  \in S_1,  \quad \sigma_2(d_i)  \in S_2.
\end{eqnarray}
Moreover, we set
\begin{eqnarray}
D'[\sigma_1, d'_i] = p_1^{-1}(\sigma_1(d'_i)) \subset
W_{|\sigma_1(\P1)} ( \simeq S_2) \subset W \\
D[\sigma_2, d_i] = p_2^{-1}(\sigma_2(d_i)) \subset
W_{|\sigma_2(\P1)}
(\simeq S_1) \subset W
\end{eqnarray}
Note that $D'[\sigma_1, d'_i]  $ and $D[\sigma_2, d_i]$ are
reduced  rational curves each of which has
one node as its singularities.
{From} (\ref{eq:p-degree1}), (\ref{eq:p-degree2})  it is easy to see
that
\begin{eqnarray}
([F], D[\sigma_2, d_i])_W = 0, \quad  ([L_i], D[\sigma_2, d_i])_W =
1,
([M_j],  D[\sigma_2, d_i])_W = 0,  \label{eq:int-l}\\
([F], D'[\sigma_1, d'_i])_W = 0, \quad  ([L_i], D'[\sigma_1,
d'_i])_W = 0,
([M_j],  D'[\sigma_1, d'_i])_W = 1. \label{eq:int-m}
\end{eqnarray}

We denote by  $\cI(k_i, \sigma_2(d_i))$ an ideal sheaf on
$S_2$  such that the quotient sheaf
$$
\cO_{S_2}/\cI(k_i, \sigma_2(d_i))
$$ is  supported on the point $\sigma_2(d_i)$ and
$ \length \cO_{S_2}/ \cI(k_i, \sigma_2(d_i))= k_i$.
We call such an ideal $\cI(k_i, \sigma_2(d_i))$ a punctual
ideal of colength $k_i$ supported on $\sigma_2(d_i)$.
And similarly for  $\cI( k'_j, \sigma_1(d'_j))$.
For each of $1\leq i\leq 12$ (resp. $1 \leq j \leq 12$), let
$\cI(k_i, \sigma_2(d_i))$ (resp. $\cI( k'_j, \sigma_1(d'_j))$)
be a punctual ideal of colength $k_i$ (resp. $k'_j$) supported on
$\sigma_2(d_i)$ (resp. $\sigma_1(d'_j)$).  We denote
by
$$
D[\cI(k_i, \sigma_2(d_i))] \quad (\mbox{resp.} 
D'[\cI( k'_j,\sigma_1(d'_j))] )
$$
the one-dimensional subscheme of $W$
defined by the pullback of the ideal sheaf $\cI(k_i,
\sigma_2(d_i))$
(resp. $\cI( k'_j, \sigma_1(d'_j))$) via $p_2$ (resp.
$p_1$).   Note that 
$$
D[\cI(k_i, \sigma_2(d_i))]_{red} = D[\sigma_2, d_i], \quad D'[\cI( k'_j,\sigma_1(d'_j))]_{red} = D'[\sigma_1, d'_j].
$$ 

Now we take the following subscheme of $W$:
\begin{eqnarray}
C = \sigma(\P1) +  \sum_{i = 1}^{12} 
D[\cI(k_i,\sigma_2(d_i))]  +
\sum_{j = 1}^{12} D'[\cI( k'_j, \sigma_1(d'_j))].
\label{eq:p-general}
\end{eqnarray}
This one dimensional  subscheme $C$ in (\ref{eq:p-general}) is
actually
a pseudo-section. }
\end{Example}

\begin{Definition}{\rm The pseudo-section $C$ in (\ref{eq:p-general})
is called
{\em of type}
$$
(\sigma, k_1, \cdots, k_{12}, k'_1, \cdots,  k'_{12}) \in MW(W)
\times
(\Z_{+})^{24}.
$$
}
\end{Definition}

\begin{Proposition}
Every  pseudo-section $C$ of $h:W \lra \P1$  can be written as in
$(\ref{eq:p-general})$.
\end{Proposition}

{\it Proof.}  Since $ ([F], C)_W = (F_1, (p_1)_* C))_{S_1} =  (F_2,
(p_2)_*
C))_{S_2} = 1$, it
is easy to see that there are sections $\sigma_i \in MW(S_i)$,
$$
(p_1)_* (C) = \sigma_1(\P1)  +  \mbox{fibers} , \quad (p_2)_* (C) =
\sigma_2(\P1) +
\mbox{fibers}.
$$
Then by definition of a pseudo-section, we can easily see that $C$
can be
written
in the form of (\ref{eq:p-general}) where $\sigma $ corresponds to
$(\sigma_1,
\sigma_2)$.
{}\hfill$\Box$

\vskp
Fix a type
$\mu  = (\sigma, k_1, \cdots, k_{12}, k'_1, \cdots,  k'_{12}) \in
MW(W) \times
(\Z_{+})^{24}$ of
a pseudo-section of $h:W \lra \P1$.
We would like to count the
``number''$n(\mu)$ of rational curves
which gives  the correct  contribution to the Gromov-Witten
invariant in
the formula (\ref{eq:g-w.inv}).
Since a pseudo-section of type $\mu$  is a non-reduced subscheme of
$W$
if some  $k_i$ or $k'_j$ is
greater than 1,
it is not easy to determine $n(\mu)$.  Of course, the
Gromov-Witten invariant should be defined as the number of
 $J$-holomorphic curves with a fixed homology class
after perturbing the complex structure of $W$ to a generic almost
complex structure $J$ (\cite{Mor-Math}, Theorem 3.3).  However at
this moment
 we do not know how to perturb the almost complex structure and how
a pseudo-section  $C$ of type $\mu$ arises as a limit of
 $J$-holomorphic curves.
 (Different $J$-holomorphic curves for generic $J$ may have the
 same limit in our complex structure of Schoen's Calabi-Yau 3-fold
$W$.)

 Here we propose the following:

 \begin{Conjecture}\label{c:cont}
 The contribution  $n(\mu) $ of all  pseudo-sections of type
 $\mu$ is given by
\begin{eqnarray}
n(\mu) = e(\Hilb^{\mu}_W) = \mbox{ Topological Euler number of }(
\Hilb^{\mu}_W),
\label{eq:cont}
\end{eqnarray}
where $\Hilb^{\mu}_W $ is the Hilbert scheme parameterizing
pseudo-sections $C
\subset W$
of type $\mu$.
\end{Conjecture}

Let $\C^2$ be the complex affine space of dimension 2 and denote by
$\Hilb^{k}(\C^2, 0)$ the Hilbert scheme
parameterizing the punctual ideal sheaves
$\cI \subset \cO_{\C^2} $ of colength $k$
supported on the origin $0 \in \C^2$

\begin{Lemma} \label{l:hilb-1}
Fix a type  $\mu =  (\sigma, k_1, \cdots, k_{12}, k'_1, \cdots,
k'_{12})$ of
pseudo-section. Then we  have a natural isomorphism of schemes
\begin{equation}
\Hilb^{\mu}_W \simeq  \prod_{i=1}^{12} (\Hilb^{k_i}(\C^2, 0)) \times
\prod_{j=1}^{12} (\Hilb^{k'_j}(\C^2, 0))
\label{eq:hilb-isom}
\end{equation}
\end{Lemma}

{\it Proof.}  {From} the definition of  pseudo-section $C$ of type
$\mu$ in
(\ref{eq:p-general}),
we have the natural morphism  $\varphi$ from $\Hilb^{\mu}_W$ to
\begin{eqnarray*}
\Hilb (\sigma(\P1) \subset W) \times
 \prod_{i=1}^{12} (\Hilb^{k_i}(S_2, \sigma_2(d_i))) \times
 \prod_{j=1}^{12} (\Hilb^{k'_j}(S_1, \sigma_1(d'_j))).
\end{eqnarray*}
of $W$
defined by
\begin{eqnarray*}
 \varphi(C) & = & \varphi(\sigma(\P1) +  \sum_{i = 1}^{12}
D[\cI(k_i, \sigma_2(d_i))]  +
\sum_{j = 1}^{12} D'[\cI( k'_j, \beta[\sigma_1, d'_j])])  \\
 & = &  (\sigma(\P1),   \ \cI[k_i, \sigma_2(d_i)],  \
\cI[k'_j,\alpha[\sigma_1,
 d'_j]).
\end{eqnarray*}
(Here   $\Hilb (\sigma(\P1) \subset W)$ denotes the connected
component of the
Hilbert scheme which contains the subscheme $\sigma(\P1) $ of $W$. )

Noting that $C$ is connected and $\sigma(\P1) \subset W$ has no
deformation (in
particular $ \Hilb(\sigma(\P1) \subset W) = 1 pt$) , we
can easily see that $\varphi$ is an isomorphism and obtain
(\ref{eq:hilb-isom}).
{}\hfill$\Box$

 \vskp
The following important lemma is  a  kind suggestion of Kota
Yoshioka.

\begin{Lemma} \label{l:hilb-2}
The Hilbert scheme $\Hilb^{k}(\C^2, 0)$ is irreducible scheme of
dimension $k-1$ and
$$
e(\Hilb^{k}(\C^2, 0)) = p(k)
$$
where $p(k)$ denotes the number of partitions of $k$.
\end{Lemma}

{\it Proof.}  The irreduciblity of $\Hilb^{k}(\C^2, 0) $ is
due to Brian\c{c}on \cite{B}.
Moreover $\Hilb^{k}(\C \{ x, y \}) $
has a natural $S^1$-action induced by $(x, y) \rightarrow
(t^a \cdot x, t^b \cdot y)$ for any weight $(a, b)$.
 Then for a general
choice of a weight $(a, b)$  its fixed points set becomes
just the set of monomial ideals of length $k$.
Now a standard argument shows that  the
 topological Euler number of
$\Hilb^{k}(\C^2, 0)$ is equal to the number of fixed points,
and it is an easy exercise to see that the
number of monomial ideals of $\C[x, y]$ with  colength $k$  is equal
to $p(k)$.{}\hfill$\Box$

\vskp
{From} Lemma~\ref{l:hilb-1} and Lemma~\ref{l:hilb-2}, we obtain the
following result.
\begin{Corollary}
Let $\mu$ and $n(\mu)$ as in Conjecture~(\ref{c:cont}), then we have
\begin{eqnarray}
n(\mu) = e(\Hilb^{\mu}_W) = (\prod_{i=1}^{12} p(k_i)) \cdot
(\prod_{j=1}^{12}
p(k'_j))
\label{eq:cont-result}
\end{eqnarray}
{}\hfill$\Box$
\end{Corollary}

\vskp

\begin{Lemma}\label{l:homology}
Let $\mu = (\sigma, k_1, \cdots, k_{12}, k'_1, \cdots, k'_{12}),
\mu' =
(\sigma', l_1, \cdots, l_{12}, l'_1, \cdots, l'_{12}) $ be two types
of
pseudo-sections.
Then a pseudo-section $C$ of type $\mu$ and $C'$ of type $\mu'$ have
the same homology class in $H_2(W)$ if and only if
\begin{eqnarray}
\sigma = \sigma',  \quad \sum_{i =1}^{12} k_i = \sum_{i =1}^{12}
l_i, \quad
\sum_{j =1}^{12} k'_j = \sum_{j =1}^{12} l'_j.
\label{eq:iff}
\end{eqnarray}
\end{Lemma}

\noindent
{\it Proof.} The ``only if" part is obvious.  Lemma~\ref{l:j-inj}
shows that
the first equality in (\ref{eq:iff}) is  necessary.  Moreover
noting  that  $[D[\cI(k_i, \sigma_2, d_i)]]$ is homologous to $k_i
\cdot (p_2)^{-1}(\sigma_2(d_i)) $ and
$[D'(\cI[k'_j, \sigma_1, d'_j)]] $ is homologous to $ k'_j \cdot
[(p_1)^{-1}(\sigma_1(d'_j))]
$
we have the other  implication.
{}\hfill$\Box$

\vskp
Let $\eta \in H_2(W)$ be such that $([F], \eta)_W = 1$, then in order
to
have non-vanishing contribution $n(\eta)$, $\eta$ must be the class of a
pseudo-section, so write
$\eta$ as $(\sigma, m \cdot p_2^{-1}( 1 pt), n p_1^{-1} ( 1 pt))$.
We call  $\eta$ {\em of type $(\sigma, m, n)$}.

\begin{Proposition}
For $\eta \in H_2(W)$ with $([F], \eta)_W = 1$  of type $(\sigma, m,
n)$, we have
\begin{eqnarray}
n(\eta)= n(\sigma, m, n)
:=  (\sum_{ \  k_1 + \cdots + k_{12} = n} \quad  \prod_{i =
1}^{12}p(k_i ))
(\sum_{ \ k'_1  + \cdots +k'_{12} = m} \quad  \prod_{j =
1}^{12}p(k'_j)).
\label{eq:coef}
\end{eqnarray}
Here  $k_i$ and  $k'_j$ run over the non-negative integers.

\end{Proposition}

\vskp
\noindent
{\it Proof.}  {From} Lemma~\ref{l:homology} and the remark above, we
have
\begin{eqnarray*}
n(\eta) =  \sum_{\mu} n(\mu)
\end{eqnarray*}
where the summation is taken over the types
$\mu = (\sigma, k_1, \cdots, k_{12}, k'_1, \cdots, k'_{12})$   of
pseudo-sections such that
$$
n = \sum_{i=1}^{12} k_i, \quad m = \sum_{j =1}^{12} k'_j.
$$
Combining this with (\ref{eq:cont-result}), we obtain the
assertion~(\ref{eq:coef}). {}\hfill$\Box$

\vskp
\vskp
\noindent
{\bf Proof of Theorem~\ref{t:1-s-prep}}
\vskp

Let us fix a homology class $\eta \in H_2(W)$ with $([F], \eta) =1$
of
type $(\sigma, n, m)$. From (\ref{eq:int-l}), (\ref{eq:int-m})  it
is easy to see that
\begin{eqnarray}
([F], \eta)_W  & =  &1,  \\
([L_i], \eta)_W & =  &
([L_i], j(\sigma))_W + n  \\
([M_j], \eta)_W & = & ([M_j], j(\sigma) )_W +  m.
\end{eqnarray}

We introduce parameters $z_0 = \log q_0, y_0 = \log s_0$.
Note that we have set in Theorem~\ref{t:mw-prep}
$$
\tau_1 = \sum_{j=0}^{8} z_j, \quad \tau_2 = \sum_{j=0}^8 y_j.
$$
Moreover just for notation in the proof, we set $v_l = \exp(2 \pi i
\tau_l)$
for
$l =1, 2$.
Recalling the definition of $T^{\sigma}$ (cf. (\ref{eq:tsigma})), we
have
\begin{eqnarray}
T^{\eta} & = &  \exp( 2 \pi i (t_0 [F] + \sum_{i=0}^8 z_i [L_i] +
\sum_{j=1}^8 y_j [M_j], \eta)_W)  \nonumber \\
   & = & p \cdot \exp(2 \pi i (\sum_{i=0}^8 z_i
[([L_i], j(\sigma))_W+ n] +  \sum_{i=0}^8 y_j [([M_j], j(\sigma))_W
+ m])) \nonumber \\
& = &  T^{\sigma} \cdot (v_1)^n \cdot
(v_2)^m
\end{eqnarray}
Then we have
\begin{eqnarray}
\lefteqn{\Psi_{A, 1}(p, q_0, \cdots q_8, s_0, \cdots, s_8)}
\nonumber  \\
&= & \sum_{\eta \in H_2(W), ([F], \eta) = 1} n(\eta) \Li_3(T^{\eta})
\nonumber
\\
&= &  \sum_{\sigma \in MW(W), n \geq 0, m \geq 0}  n(\sigma, n, m)
\cdot
 \Li_3( T^{\sigma} \cdot (v_1)^n \cdot (v_2)^m) \nonumber \\
&= & \sum_{N = 1}^{\infty}
[\sum_{\sigma \in MW(W), n \geq 0, m \geq 0} n(\sigma, n, m) \cdot
\frac{(T^{\sigma})^N \cdot (v_1)^{Nn}\cdot  (v_2)^{Nm}}{N^3}]
\nonumber \\
& = & \sum_{N = 1}^{\infty} \frac{1}{N^3} \cdot
[\sum_{\sigma \in MW(W)} (T^{\sigma})^N] \cdot
[ \sum_{n \geq 0, m \geq 0} n(\sigma, n, m)
(v_1)^{Nn} \cdot (v_2)^{Nm}].  \nonumber \\
  &  &   \label{eq:final}
\end{eqnarray}
Note that the last equality follows from the fact that
$n(\sigma, n, m)$ does not depend on $\sigma$.

On the other hand, from  equality~(\ref{eq:coef}) we have
$$
[(\sum_{k =0}^{\infty} p(k) (v_1)^k )(\sum_{k' =0}^{\infty} p(k')
(v_2)^{k'})]^{12}
= \sum_{n \geq 0, m \geq 0} n(\sigma, n, m) (v_1)^n \cdot (v_2)^m.
$$
Moreover as in the proof of Theorem~\ref{t:mw-prep} we can see that
$$
[\sum_{\sigma \in MW(W)} (T^{\sigma})^N] =
( p \prod_{i=1}^8 (q_i \cdot s_i))^N
\Theta_{E_8}^{root}(N  \tau_1, N\cdot \z) \cdot
\Theta_{E_8}^{root}(N \tau_2, N \cdot \y).
$$
Combining these equalities with (\ref{eq:final}), we
obtain the proof of Theorem~\ref{t:1-s-prep}. {}\hfill$\Box$

\vskp

\section{The restricted A-model Prepotential}
\label{s:rest}

In order to compare the prepotential of the A-model Yukawa coupling
of $W$
with the B-model Yukawa coupling
of the mirror partner $W^*$,
which we obtain in Section~\ref{s:B-model},
we need to take a special restriction of the variables of the
prepotential,  that is, we have to specify
the parameters which correspond to the line bundles
which are induced from the ambient space $\P1 \times \P2 \times
\P2$.
Let $\iota:W \hookrightarrow
\P1 \times \P2 \times \P2$ be the natural embedding.
Then we set
\begin{equation}
[F] = \pi_1^*(\cO_{\P1}(1)), \quad [H_1]=  \pi_2^*(\cO_{\P2}(1)),
\quad [H_2]
=  \pi_3^*(\cO_{\P2}(1)), \label{eq:res-picard}
\end{equation}
and introduce corresponding parameters as follows:
\begin{eqnarray}
  \quad  [F]   \leftrightarrow   p = U_0 = \exp (2 \pi i t_0),
\nonumber \\
 \quad [H_1]  \leftrightarrow  U_1 = \exp (2 \pi i t_1),
\label{eq:sp-coord}
 \\
 \quad  [H_2]  \leftrightarrow  U_2 = \exp (2 \pi i t_2).  \nonumber
\end{eqnarray}
Now we consider the following restricted prepotential
\begin{eqnarray}
\Psi_{A}^{res} = \mbox{topological term} + \sum_{0 \neq \eta \in
H_2(W)}
n(\eta) \Li_3(U^{\eta})
\label{eq:res-fullprep}
\end{eqnarray}
where
\begin{eqnarray}
U^{\eta} & = &  \exp( 2 \pi i (t_0 [F] + t_1 [H_1] + t_2 [H_2],
\eta)_W) \\
   & = & p^{([F], \eta )_W } \cdot (U_1)^{([H_1], \eta )_W } \cdot
(U_2)^{ ([H_2], \eta )_W}.
\end{eqnarray}

Moreover, we can define the $k$-sectional part and the
Mordell-Weil part of the
restricted prepotential
by
\begin{eqnarray}
\Psi_{A, k}^{res} = \sum_{0 \not\equiv \eta
\in H_2(X, \Z), \ (F, \eta)= k}  n(\eta) \Li_3( U^{\eta}),
\label{eq:res-sec-prep}
\end{eqnarray}
\begin{eqnarray}
\Psi_{A, MW(W)}^{res} = \sum_{\sigma \in MW(W)}
\Li_3(U^{j(\sigma)}),
\label{eq:res-prep}
\end{eqnarray}
respectively.

\begin{Proposition}
\begin{equation}
\Psi_{A, MW(W)}^{res}(p, t_1, t_2)
= \sum_{n=1}^{\infty} \frac{p^n}{n^3} \cdot
\Theta_{E_8}(3n t_1,  n t_1 \gamma) \cdot \Theta_{E_8}(3n
t_2,  n t_2 \gamma) \label{eq:resprep}
\end{equation}
\begin{equation}
\Psi_{A, 1}^{res}(p, t_1, t_2)
= \sum_{n=1}^{\infty} \frac{p^n}{n^3} \cdot
A^{res}(n t_1) \cdot A^{res}(n t_2)
\label{eq:A-res-prep}
\end{equation}
where
$$
\gamma = (1, 1, 1, 1, 1, 1, 1, -1)
$$
and
\begin{eqnarray}
A^{res}(t) &=&  \Theta_{E_8}( 3 t, t \cdot \gamma) \cdot
(\sum_{n=0}^{\infty}
p(n)\exp (2\pi i n(3 t))^{12}  \\
& = & \Theta_{E_8}( 3 t, t \cdot \gamma) \cdot
\frac{\exp( 3 \pi i t)}{[\eta( 3t )]^{12}}   \nonumber \\
& = & \Theta_{E_8}( 3 t, t \cdot \gamma) \cdot
\frac{1}{[\prod_{m \geq 1 }(1 - \exp(2 \pi i m (3t)))]^{12}}.
\nonumber
\label{eq:ares}
\end{eqnarray}
\label{prop:res}
\end{Proposition}

\noindent
{\it Proof.}   {From}  Relation~(\ref{eq:h-f}), we obtain for every
$\sigma \in
MW(W)$
\begin{eqnarray*}
\lefteqn{([H_1], j(\sigma))_W  = (H_1, j(\sigma_1))_{S_1} }\\
   & &  =(2 F_1 + 5 L_0 -2 L_1 - L_2   + L_8, j(\sigma_1))_{S_1} \\
   &  & = 2 + 1/2 (5 Q_1(\sigma_1) - 2 Q_1(\sigma_1 + a_1) -
Q_1(\sigma_1 +
   a_2) +  Q_1(\sigma_1 + a_8)) \\
  & &  = \frac{3}{2} Q_1(\sigma_1) + B_1(\sigma_1, -2 a_1 -a_2 +
a_8),
\end{eqnarray*}
and a similar equation for $([H_2], j(\sigma))_W$.
Then from Remark~\ref{rem:root}, we see that
\begin{eqnarray*}
\gamma = -2 a_1 - a_2 + a_8 & = 
&  -[(\e_1 + \e_8) -(\e_2 + \e_3 + \e_4 +
\e_5 + \e_6
+ \e_7)] + \\
  & & -(\e_2 - \e_1) + (\e_1 + \e_2) \\
  &=& \e_1 + \e_2 + \e_3 + \e_4 
+ \e_5 + \e_6 + \e_7 - \e_8
\end{eqnarray*}
Therefore we see that
\begin{eqnarray*}
\lefteqn{U^{j(\sigma)}  = \exp (2 \pi i (t_0 + ([H_1],j(\sigma))_W
t_1 +
([H_2],j(\sigma))_W t_2))} \\
 &&= p \cdot \exp (2 \pi i(\frac{3t_1}{2} Q_1(\sigma_1)+B_1(
\sigma_1, t_1
 \gamma ) ) \exp ( 2 \pi i( \frac{3t_2}{2}Q_2(\sigma_2)+B_2(
\sigma_2, t_2
 \gamma ) ) )
\end{eqnarray*}
Then as in the proof of Theorem~\ref{t:mw-prep}, we can obtain
Assertion~(\ref{eq:resprep}).  The proof of
Assertion~(\ref{eq:A-res-prep}) is
similar.
{}\hfill$\Box$

\vskp
Now we consider the expansion of (\ref{eq:resprep}) and
(\ref{eq:A-res-prep}) with respect to the variable $p$. 
\begin{equation}
\Psi_{A, MW(W)}^{res} = 
 p \cdot \Theta_{E_8}(3 t_1,   t_1 \gamma)
\cdot
\Theta_{E_8}(3 t_2,   t_2 \gamma) + O(p^2)
\label{eq:2ptheta}
\end{equation}
\begin{equation}
\Psi_A - \mbox{topological term} =
 \Psi_{A, 1}^{res} + O(p^2) =  
p \cdot A^{res} (t_1) \cdot
A^{res}(t_2)+ O(p^2)
\label{eq:2-A-prep}
\end{equation}

\nvskp
Let us define sequences of positive integers $\{c_n \} $ and
$\{ a_n \}$ by
\begin{eqnarray}
\Theta_{E_8}(3 t,   t \gamma) = \sum_{n = 0}^{\infty} c_m \exp (2
\pi m i t) =
\sum_{n = 0}^{\infty} c_m U^m
\label{eq:expand} \\
A^{res}( t ) = \sum_{n = 0}^{\infty} a_m \exp (2 \pi m i t) =
\sum_{n =
0}^{\infty} a_m U^m.
\label{eq:expand-a}
\end{eqnarray}

Let us also expand  functions
as follows:
\begin{eqnarray}
p \cdot \Theta_{E_8}(3 t_1,   t_1 \gamma) \cdot \Theta_{E_8}(3 t_2,
 t_2
\gamma) =
\sum_{n_1 \geq 0, n_2 \geq 0} M_{1, n_1,n_2} \cdot p \cdot
(U_1)^{n_1} (U_2)^{n_2}\label{eq:expandprod1} \\
p \cdot A^{res}(t_1) \cdot A^{res}(t_2) = \sum_{n_1 \geq 0, n_2 \geq 0}
N_{1,
n_1,n_2} \cdot p \cdot
(U_1)^{n_1} (U_2)^{n_2}.
 \label{eq:expandprod2}
\end{eqnarray}

The proof of the following proposition follows from
(\ref{eq:expandprod1}),
(\ref{eq:expandprod2})
and Proposition~\ref{prop:res}.

\begin{Proposition}
\begin{enumerate}
\item
\begin{eqnarray}
 M_{1,n_1,n_2} =  c_{n_1} \cdot c_{n_2}, \quad
N_{1, n_1, n_2} = a_{n_1} \cdot a_{n_2}
\label{eq:split}
\end{eqnarray}

\item Let $f:S \lra \P1$ be  a generic rational elliptic surface as
in
section~\ref{s:schoen}. Then
in the expansion~$(\ref{eq:expand})$ the coefficient $c_m$ is the
number of
sections of $f:S \lra \P1$ with degree $(H, [\sigma(\P1)])_S = m$,
that is,
$$
c_m = \sharp \{ \ \ \sigma \in MW(S) \quad | \quad  (H, [\sigma(\P1)])_S = m
  \ \ \},
$$
where $H$ is the class of the  total transform of a line on $\P2$.

\item  The integer  $M_{1,n_1, n_2}$ $($resp.$N_{1,n_1, n_2}$$)$
is the number of sections $\sigma \in MW(W)$  (resp. pseudo-sections
$\eta$)
of $h:W \lra \P1$ with bidegree $(n_1, n_2) $  where $n_i = ([H_i],
j(\sigma))_W$
$($resp. $n_i = ([H_i], \eta)_W$ $)$ is
the degree with respect to $[H_i]$.

\end{enumerate}
\end{Proposition}
{}\hfill$\Box$

\begin{Remark}{\rm
The factorization property of  $M_{1,n_1, n_2}$ and $N_{1, n_1,
n_2}$
in $(\ref{eq:split})$  follows from the fact that sections and
pseudo-sections
of $h:W \lra \P1$ can be split as in  $(\ref{eq:isom})$,
(\ref{eq:p-general}).
}
\end{Remark}
\begin{Remark}{\rm
Note that the sequences $\{ c_m \}$ and $\{a_m \}$ are connected to
each other
by the formula:
$$
\sum_{n = 0}^{\infty} a_m U^m = [\sum_{n = 0}^{\infty} c_m U^m][
\sum_{k
=0}^{\infty} p(k) U^{3k}]^{12}.
$$
The number $a_m$ can be considered as the number of pseudo-sections
$C$
of a
generic rational elliptic surface $f:S \lra \P1$ of degree $m$ with
respect to
the divisor class $[H]$.
The term
$$
[ \sum_{k =0}^{\infty} p(k) U^{3k}]^{12}
$$
is nothing but the contribution of 12 singular fibers of type $I_1$,
when we
count the contribution of one singular fiber of type $I_1$ with
multiplicity
$k$ as $p(k)$.
}
\end{Remark}

Here we will expand $\Theta_{E_8}(3 t,   t \gamma)$ and $A^{res}(t)$
in the variable $U = \exp(2 \pi i t)$ and give the table of
coefficients $c_n$ and $a_n$ up to order 50. (See also the last
remark  of Section~\ref{s:B-model}).   We can use
Proposition~\ref{pr:e8theta} to obtain the following expansion.

\pagebreak
\begin{center}
Table. 1
\end{center}
$$
\Theta_{E_8} (3 t ,   t \gamma )  =  \sum_{m \geq 0} c_m U^m.
$$
$$
\begin{array}{lclclcl}
\begin{array}{r|l}
n&c_n\\
&\\
0& 9 \\
1& 36\\
2 & 126\\
3 & 252 \\
4 & 513 \\
5 & 756 \\
6 & 1332\\
7 & 1764 \\
8 &  2808 \\
9& 3276 \\
10 &4914 \\
11 &  5616\\
12 & 8190 \\
\end{array}
&\hspace{5mm}&
\begin{array}{r|l}
n&c_n\\
&\\
13 &  8892 \\
14 &  12168 \\
15 &  13104\\
16 & 17766 \\
17 &  18648 \\
18 &  24390 \\
19 &  25200\\
20 & 33345  \\
21 &  33516\\
22 &  43344 \\
23 &  43092  \\
24 &  55692 \\
25 & 54684\\
\end{array}
&\hspace{5mm}&
\begin{array}{r|l}
n&c_n\\
&\\
26 & 68922 \\
27 & 68796   \\
28 &  86580  \\
29 &  84168 \\
30 &  103824 \\
31 & 101556    \\
32 &  127647  \\
33 &  121212    \\
34 &  148878   \\
35 &  143964    \\
36 &  178776     \\
37 & 170352    \\
38 &   205380    \\
\end{array}
&\hspace{5mm}&
\begin{array}{r|l}
n&c_n\\
&\\
39 &   197136    \\
40 &  241920\\
41 &  227556   \\
42 &  276948  \\
43 & 262080   \\
44 & 319410    \\
45 &  298116  \\
46 &  357912\\
47 &  341460   \\
48 &  410958   \\
49 &  382356   \\
50 & 458208  \\
  & \\
\end{array}
\end{array}
$$

\begin{center}
{\bf Table 2. (Table for $\{ a_n \} $).}
\end{center}
$$
A^{res}(t) =  \sum_{m \geq 0} a_m U^m = (\sum_{n \geq 0} c_n U^n)
\cdot
(\sum_{k \geq 0} p(k) U^{3k} )^{12}.
$$
$$
\begin{array}{lclcl}
\begin{array}{r|l}
n&a_n\\
&\\
0&9\\
 1& 36\\
 2 & 126\\
 3 & 360\\
 4 & 945\\
 5 & 2268\\
6 & 5166\\
7 & 11160\\
8 & 23220\\
9 & 46620\\
10 & 90972\\
11 & 172872\\
12 & 321237 \\
13 & 584640\\
14 & 1044810\\
15 & 1835856\\
16 & 3177153\\
\end{array}
&\hspace{5mm}&
\begin{array}{r|l}
n&a_n\\
&\\
17 & 5421132\\
18 & 9131220\\
19 & 15195600\\
20 & 25006653\\
21 & 40722840\\
22 & 65670768\\
23 & 104930280\\
24 & 166214205\\
25 & 261141300\\
26 & 407118726\\
27 & 630048384\\
28 & 968272605\\
29 & 1478208420\\
30 & 2242463580\\
31 & 3381344280\\
32 & 5069259342\\
33 & 7557818940\\
\end{array}
&\hspace{5mm}&
\begin{array}{r|l}
n&a_n\\
&\\
34 & 11208455370\\
35 & 16538048640\\
36 & 24282822798\\
37 & 35487134928\\
38 & 51626878470\\
39 & 74779896240\\
40 & 107861179482\\
41 & 154945739844\\
42 & 221711362038\\
43 & 316042958880\\
44 & 448856366490\\
45 & 635216766732\\
46 & 895854679650\\
47 & 1259213600736\\
48 & 1764210946995\\
49 & 2463949037340\\
50 & 3430694064888
\end{array}
\end{array}
$$
\vskp

\pagebreak

\section{The prepotential of the B-model Yukawa coupling}
\label{s:B-model}

In this section we study the prepotential of the B-model Yukawa
coupling
for the mirror $W^*$ of Schoen's example
in the sense of Batyrev-Borisov~\cite{Batyrev-Borisov}
 and compare it with the prepotential
for the A-model Yukawa coupling of $W$.
Formula (\ref{eq:Psi}) gives this B-model prepotential $\Psi_{B}$
explicitly.

In order to determine the B-model prepotential
we will basically follow the recipe of 
\cite{HKTY,Sti} which uses only
the
toric data of {the} A-model side. 
However in order to give an
intuitive picture of the mirror $W^*$ we will put
here the  orbifold  construction of the 
mirror $W^*$ of Schoen's example
and the Picard-Fuchs equations of the periods of a
holomorphic 3-form of $W^*$.

Based on the Batyrev-Borisov mirror construction (cf.
\cite{Batyrev-Borisov}, \cite{HKTY}) for complete intersection
Calabi-Yau manifolds in toric varieties
we can derive the following

\begin{Proposition}\label{orbifold}
The family of
mirror Calabi-Yau 3-folds of $W$ is obtained by the orbifold
construction with group  $\Z_3\times\Z_3$ from
the subfamily $W_{\alpha_0, \alpha_1, \beta_0, \beta_1}$ of $W$:
\begin{eqnarray*}
\lefteqn{W_{\alpha_0, \alpha_1, \beta_0, \beta_1} = } \\
&& \left\{  \;
[z_0:z_1] \times [ x_0: x_1: x_2]\times[ y_0: y_1: y_2]
 \in  \P1 \times \P2 \times \P2 
\;|\;P_1 = P_2 = 0 \;
\right\}
\end{eqnarray*}
where
$$
\begin{array}{rcl}
P_1&=&(x_0^3+x_1^3+x_2^3+\alpha_0x_0x_1x_2)z_1+\alpha_1x_0x_1x_2 z_0
\;,\\
P_2&=&(y_0^3+y_1^3+y_2^3+\beta_0y_0y_1y_2)z_0+\beta_1y_0y_1y_2 z_1
\;  \\
\end{array}
$$
and the group $\Z_3\times\Z_3$ is generated by
\begin{eqnarray} \label{eq:action}
\lefteqn{g_1: ([z_0:z_1],[x_0:x_1:x_2],[y_0:y_1:y_2]) \hspace{2cm}
}  \nonumber \\
    & & \hspace{1cm} \mapsto([z_0:z_1],[x_0:\omega x_1:\omega^2
x_2],[y_0:y_1:y_2]),  \\
\lefteqn{ g_2: ([z_0:z_1],[x_0:x_1:x_2],[y_0:y_1:y_2]) \hspace{2cm}
} \nonumber  \\
       && \hspace{1cm} \mapsto ([z_0:z_1],[x_0:x_1:x_2],[y_0:\omega
y_1:\omega^2 y_2]), \nonumber
\end{eqnarray}
with $\omega=e^{2\pi i /3}$. That is, the mirror $W^*$ is
$$
W^* = W_{\alpha_0, \alpha_1, \beta_0, \beta_1}/(\Z_3 \times \Z_3).
$$
\end{Proposition}
\noindent
{\it Proof.} See Section \ref{s:B-model equations}, Appendix I.
{}\hfill$\Box$

In the equations $P_1$ and $P_2$ above we have kept four parameters
$\alpha_0,\alpha_1, \beta_0,\beta_1$ for symmetry reasons. However
only three of them are essential because of the scaling of the
variables $z_0,z_1$. After the orbifoldization this three parameter
deformation describes a three dimensional subspace in the complex
structure (B-model) moduli space of $W^*$. The full complex
structure moduli
space has dimension 19. Under the mirror symmetry the three
dimensional subspace will be mapped to the subspace in the
complexified K\"ahler moduli space parameterized by $(t_0, t_1,t_2)$
in (\ref{eq:sp-coord}).
The B-model calculations are local calculations based on the
variation of
the Hodge structure for the family $W^*$ about the {\it Large
complex
structure limit} (LCSL). A mathematical characterization of LCSL
is
given in \cite{Mor-II}. Here we simply follow a general recipe
applicable to CICYs in toric varieties to find a LCSL and write
the Picard-Fuchs differential equations \cite{HKTY}. We find that
the
origin
of the local coordinate system $u=(u_0,u_1,u_2)$ with
$u_0={\alpha_1\beta_1 \over \alpha_0 \beta_0}, \;
u_1=-{1\over \alpha_0^3}$ and $u_2=-{1\over \beta_0^3}$ is a LCSL,
and that the Picard-Fuchs (PF) differential operators about this
point are
\begin{equation}
\begin{array}{rcl}
D_1 &= & (3\theta_{u_1}-\theta_{u_0})\theta_{u_1}
   -9u_1(3\theta_{u_1}+\theta_{u_0}+2)(3\theta_{u_1}+\theta_{u_0}+1)
\\
  & & + u_0\theta_{u_1}(3\theta_{u_2}+\theta_{u_0}+1) \;\;, \\[.5em]
D_2&=& (3\theta_{u_2}-\theta_{u_0})\theta_{u_2}
   -9u_2(3\theta_{u_2}+\theta_{u_0}+2)(3\theta_{u_2}+\theta_{u_0}+1)
\\
& &   + u_0\theta_{u_2}(3\theta_{u_1}+\theta_{u_0}+1) \;\;, \\[.5em]
D_3 &=& \theta_{u_0}^2-u_0(3\theta_{u_1}+\theta_{u_0}+1)
    (3\theta_{u_2}+\theta_{u_0}+1)     \;\;,
\end{array}
\label{eq:PF}
\end{equation}
with $\theta_{u_i}=u_i{\partial \; \over \partial u_i}$.
We note that if we set $u_0=0$ in (\ref{eq:PF}), the operators
$D_1$ and $D_2$ reduce to the PF equations for the Hesse pencil of
elliptic curves.
Local solutions about $u=0$ have several interesting properties.
To state these, we denote the three elements $[F], [H_1]$ and
$[H_2]$
in the Picard group Pic$(W)$ by $J_0, J_1$ and $J_2$, respectively.
By the notation $K_{ijk} \; (i,j,k=0,1,2)$ we denote the classical
intersection numbers among the corresponding divisors. Then the
nonzero components are calculated,  up to obvious permutations of the
indices,  by
\begin{equation}
K_{012}=9 \;\;,\;\; K_{112}=K_{122}=3 \;\;.
\label{eq:Kcl}
\end{equation}

\begin{Proposition}
\begin{enumerate}
\item
The Picard-Fuchs equation (\ref{eq:PF}) has only one regular
solution, namely
\begin{equation}
\Omega^{(0)}(u):=\sum_{m_0,m_1,m_2\geq 0}
{(m_0+3m_1)!\;(m_0+3m_2)! \over (m_0!)^2\, (m_1!)^3 \, (m_2!)^3 }
u_0^{m_0}u_1^{m_1}u_2^{m_2}
\label{eq:wo}
\end{equation}
\item
All other solutions of (\ref{eq:PF}) have logarithmic regular
singularities and
have
the following form in terms of the classical Frobenius method
\begin{equation} \label{eq:wi}
\begin{array}{rcl}
&\displaystyle{
\Omega^{(1)}_i(u):=
\drho{i}
\Omega(u,\rho)\vert_{\rho=0}} \;,\hspace{4mm}  \\
&\displaystyle{ \Omega^{(2)}_i(u):=
{1\over2}
\sum_{j,k=0,1,2}K_{ijk}\drho{j}\drho{k}\Omega(u,\rho)\vert_{\rho=0}
\;,\; }\\
&\displaystyle{
\Omega^{(3)}(u) :=
-{1\over3!}\sum_{i,j,k=0,1,2}K_{ijk}\drho{i}\drho{j}\drho{k}
\Omega(u,\rho)\vert_{\rho=0} \;,}
\end{array}
\end{equation}
with
\end{enumerate}
$$
\Omega(u,\rho):= \hspace{10.5cm}
$$
$$
\displaystyle{ \sum_{m_0,m_1,m_2\geq 0}
{ (1+\rho_0+3\rho_1)_{m_0+3m_1}  (1+\rho_0+3\rho_2)_{m_0+3m_2}
\over
(1+\rho_0)_{m_0}^2 (1+\rho_1)_{m_1}^3 (1+\rho_2)_{m_2}^3 }
 u_0^{m_0+\rho_0}u_1^{m_1+\rho_1}u_2^{m_2+\rho_2}}
$$
and $K_{ijk}$ being the coupling in $(\ref{eq:Kcl})$. The notation
$(x)_m$ represents the Pochhammer symbol:
$(x)_m:=x(x+1)\cdots(x+m-1)\,.$
{}\hfill$\Box$
\end{Proposition}

Now we are ready to define the B-model prepotential and the mirror
map:

\begin{Definition} {\rm We define the {\em B-model prepotential} by
\begin{equation}\label{eq:Psi}
\Psi_B(u)={1\over2}\left({1\over\Omega^{(0)}(u)}\right)^2
\Bigl\{ \Omega^{(0)}(u)\Omega^{(3)}(u) + \sum_i
\Omega_i^{(1)}(u)\Omega_i^{(2)}(u) \Bigr\} \;\;.
\end{equation}}
\end{Definition}

\begin{Definition} \label{d:mirrormap}{\rm
We define the {\em special coordinates on
the B-model moduli space} by
\begin{equation}\label{eq:coord}
t_j={1\over 2\pi i}{\Omega^{(1)}_j(u) \over \Omega^{(0)}(u) }
\:,\hspace{4mm}U_j:=e^{2\pi i t_j} \hspace{4mm}(j=0,1,2)\;.
\end{equation}
Then $U_0,U_1,U_2$ are functions of $u_0,u_1,u_2$ and
$U_j=u_j+\mbox{higher
order terms}$.
The inverse map $(u_0(U),\,u_1(U),\,u_2(U))$
is called {\em the mirror map}.}
\end{Definition}

\begin{Conjecture} \label{c:mirror}
 {\bf (Mirror Conjecture)}
The B-model prepotential $\Psi_B(u)$ combined with the mirror map
has the expansion
\begin{equation}\label{eq:mircon}
\Psi_B(u(U))
={(2\pi i)^3 \over 3!}\sum_{i,j,k=0,1,2} K_{ijk}t_it_jt_k
+\sum_{n_0,n_1,n_2 \geq 0} N_{n_0,n_1,n_2}
\Li_3(U_0^{n_0}U_1^{n_1}U_2^{n_2})
\end{equation}
where $N_{n_0,n_1,n_2}$ is the number of rational curves
$\varphi: \BP^1 \mapsto W$ with
$(J_i,\varphi_*([\BP^1]))$ $= n_i, \;\;(i=0,1,2)$.
In our context, we can state the conjecture in more precise form as
follows:
\begin{eqnarray}
\Psi^{res}_A(U_0, U_1, U_2) = \Psi_B(u(U_0, U_1, U_2))
\label{eq:mirror}
\end{eqnarray}
where $\Psi^{res}_A(U_0, U_1, U_2)$ is the restricted $A$-model
prepotential
defined in $(\ref{eq:res-fullprep})$.
\end{Conjecture}

\vskp
Next we briefly sketch the approach of \cite{Sti} for calculating
the
B-model prepotential by using only toric data of the A-model
side.
 This starts from the
observation that Schoen's example $W$ can be embedded in
$\P1\times \P2\times \P2$ as the intersection of a hypersurface of
degree $(1,3,0)$ and a hypersurface of degree $(1,0,3)$. (cf.
Section~\ref{s:schoen}).
So $W$ is the zero locus of a (general) section of the rank $2$
vector
bundle $\cO (1,3,0)\oplus \cO (1,0,3)$ on
$\P1\times \P2\times \P2$.
This vector bundle can be constructed as a quotient of an open part
of
$\C^{10}$ by a $3$-dimensional subtorus of $(\C^\ast)^{10}$ acting
by
coordinatewise multiplication. The subtorus is the image of the
homomorphism
$(\C^\ast)^{3}=\Z^{3}\otimes\C^\ast\rightarrow
\Z^{10}\otimes\C^\ast=(\C^\ast)^{10}$ given by the $3\times
10$-matrix
\begin{equation}\label{eq:Bmat}
\sB:=\left(
\begin{array}{rrrrrrrrrr}
-1&-1&1&1&0&0&0&0&0&0\\
-3&0&0&0&1&1&1&0&0&0\\
0&-3&0&0&0&0&0&1&1&1
\end{array}
\right)
\end{equation}
The open part of $\C^{10}$ is
\begin{equation}\label{eq:union}
\bigcup_{(i,j,k)\in\{3,4\}\times\{5,6,7\}\times\{8,9,10\}}
\;\C^{10}_{(i,j,k)}
\end{equation}
with
$$
\C^{10}_{(i,j,k)}:=\{(x_1,\ldots,x_{10})\in\C^{10}\;|\;
x_i\neq 0\,,\:x_j\neq 0\,,\:x_k\neq 0\;\}
$$
We view $\{3,4\}\times\{5,6,7\}\times\{8,9,10\}$ as a collection of
18 subsets
of $\{1,\ldots,10\}$ and note that the complement of the union
of these 18 subsets is $\{1,2\}$. As explained in \cite{Sti}
this bit of combinatorial input suffices to explicitly write down
the
hypergeometric function from which one can subsequently compute the
B-model
prepotential.

This hypergeometric function is a priori a function in 10 variables
$v_1,\ldots,v_{10}$, which correspond to the a priori 10
coefficients in the
equations $P_1$ and $P_2$ in proposition \ref{orbifold}:
\begin{eqnarray*}
\Phi&:=&
(\bar J_0+3\bar J_1)(\bar J_0 + 3\bar J_2)\times v_1^{-1}v_2^{-1}
\times\:u_1^{\bar J_0} u_2^{\bar J_1} u_3^{\bar J_2}\times\\
&\times&
\sum_{m_0,m_1,m_2\geq 0}
\frac{(1+\bar J_0+3\bar J_1)_{m_0+3m_1}\cdot
(1+\bar J_0 + 3\bar J_2)_{m_0+3m_2}}{
{(1+\bar J_0)_{m_0}}\!^2\cdot {(1+\bar J_1)_{m_1}}\!^3\cdot
{(1+\bar J_2)_{m_2}}\!^3}\:u_0^{m_0} u_1^{m_1} u_2^{m_2}
\end{eqnarray*}
with
$$
u_0 \::=\,v_1^{-1}v_2^{-1}v_3v_4\:,\hspace{3mm}
u_1 \::=-\,v_1^{-3}v_5v_6v_7\:,\hspace{3mm}
u_2 \::=-\,v_2^{-3}v_8v_9v_{10}
$$
and where $\bar J_0,\,\bar J_1,\,\bar J_2$ are elements in the ring
$$
{\cal R}_{\BP^1\times\BP^2\times\BP^2}:=
\Z[\bar J_0,\bar J_1,  \bar J_2]/
(\bar J_0^2,\bar J_1^3,\bar J_2^3).
$$
So, $v_1v_2\Phi$ is an element of
$$
\left((\bar J_0+3\bar J_1)(\bar J_0 + 3\bar J_2)
{\cal R}_{\BP^1\times\BP^2\times\BP^2}\right)\otimes
\Q[[u_1,u_2,u_3]][\log u_1,\,\log u_2,\,\log u_3]\,.
$$

The map multiplication by $(\bar J_0+3\bar J_1)(\bar J_0 + 3\bar
J_2)$
on ${\cal R}_{\BP^1\times\BP^2\times\BP^2}$ induces an isomorphism
of
linear spaces from the ring
$$
 {\cal  R}_{toric}:= {\cal R}_{\BP^1\times\BP^2\times\BP^2}/
  Ann((\bar J_0+3\bar J_1)(\bar J_0 + 3\bar J_2))
$$
onto the ideal
$(\bar J_0+3\bar J_1)(\bar J_0 + 3\bar J_2)
 {\cal R}_{\BP^1\times\BP^2\times\BP^2}$.

${\cal R}_{\BP^1\times\BP^2\times\BP^2}$ is in fact the cohomology
ring of the
ambient toric variety $\P1\times\P2\times\P2$ and ${\cal
R}_{toric}$
is a subring of the Chow ring of $W$. The classes of  $\bar J_0,\bar
J_1,\bar
J_2$ in ${\cal  R}_{toric}$ correspond to the elements $J_0,J_1,J_2$
of
Pic$(W)$ defined earlier.
One easily checks that ${\cal  R}_{toric}$ is a free $\Z$-module of
rank $8$
with basis
$\{
1\:,\;J_0\:,\;J_1\:,\;J_2\:,\;
J_1^2\:,\;J_1J_2\:,\;J_2^2\:,\;
J_1^2J_2
\}$
and that the following relations hold
\begin{eqnarray*}
&&
J_0^2=J_1^3=J_2^3=J_0J_1^2=J_0J_2^2=0\,,
\\
&&
J_0J_1=3J_1^2\:,\hspace{3mm}J_0J_2=3J_2^2\:,\hspace{3mm}
J_1J_2^2=J_1^2J_2\:,\hspace{3mm}J_0J_1J_2=3J_1^2J_2
\end{eqnarray*}
Instead of $\bar\Phi$ we may as well work with
\begin{eqnarray*}
\Omega(u,J):=  \hspace{9cm}
\end{eqnarray*}
$$
\sum_{m_0,m_1,m_2\geq 0}
\frac{(1+ J_0+3 J_1)_{m_0+3m_1}\cdot
(1+ J_0 + 3 J_2)_{m_0+3m_2}}{
{(1+ J_0)_{m_0}}\!^2\cdot {(1+ J_1)_{m_1}}\!^3\cdot
{(1+ J_2)_{m_2}}\!^3}\:u_0^{m_0+J_0} u_1^{m_1+J_1} u_2^{m_2+J_2}
$$
Using the notations (\ref{eq:Kcl}), (\ref{eq:wo}), (\ref{eq:wi})
and
$J_0^\vee:={1\over9}J_1J_2-{1\over27}J_0J_1-{1\over27}J_0J_2, \;$
$J_1^\vee:={1\over9}J_0J_2,\;$ $ J_2^\vee:={1\over9}J_0J_1$ and
$vol:={1\over9}J_0J_1J_2$
the relation between the two approaches may now be formulated as
\begin{Proposition}
$$
\Omega(u,J)=\Omega^{(0)}(u)+\sum_{i=0}^2 \Omega^{(1)}_i(u) J_i
+\sum_{i=0}^2 \Omega^{(2)}_i(u) J_i^\vee -
 \Omega^{(3)}(u) vol
$$
{}\hfill$\Box$
\end{Proposition}

$\Omega(u,J)$ is an element of the ring
${\cal  R}_{toric}\otimes
\Q[[u_1,u_2,u_3]][\log u_1,\,\log u_2,\,\log u_3]\,.$ It is $1$
modulo
$J_0,J_1,J_2,$ $u_0,u_1,u_2$ and hence its logarithm
also exists in the ring
${\cal  R}_{toric}\otimes
\Q[[u_1,u_2,u_3]][\log u_1,\,\log u_2,\,\log u_3]\,.$ Expanding
$\log\Omega(u,J)$ with respect to the basis $\{1,\,J_0,\,J_1,\,J_2,
J_0^\vee,\,J_1^\vee,\,J_2^\vee,\,vol\}$ of ${\cal  R}_{toric}$ one
finds
$$
\log\Omega(u,J)=\log\Omega^{(0)}(u)\,+\,
\sum_{j=0}^2\log U_j\,J_j\,+\,
\sum_{j=0}^2P_j\,J_j^\vee\,+\,P\,vol
$$
with $U_j$ as in (\ref{eq:coord}) and hence $\log U_j= 2\pi i t_j$.
A straightforward computation shows (see also  (\ref{eq:Psi}) and
(\ref{eq:mircon}))
\begin{eqnarray}
\nonumber
P&=&- \left(\frac{\Omega^{(3)}(u)}{\Omega^{(0)}(u)}\,+\,
\sum_{j=0}^2\frac{\Omega_j^{(1)}(u)}{\Omega^{(0)}(u)}
\frac{\Omega_j^{(2)}(u)}{\Omega^{(0)}(u)}\,-\,
\frac{(2\pi i)^3}{3}\sum_{m,j,k=0,1,2} K_{mjk}t_mt_jt_k
 \right)\\
\nonumber
&=&-{2} \left(\Psi_B(u)\,-\,
\frac{(2\pi i)^3}{3!}\sum_{m,j,k=0,1,2} K_{mjk}t_mt_jt_k
 \right)\\
\label{eq:P}
&=&-{2} \sum_{n_0,n_1,n_2 \geq 0} N_{n_0,n_1,n_2}
\Li_3(U_0^{n_0}U_1^{n_1}U_2^{n_2})
\end{eqnarray}

\begin{Proposition}
Let the numbers $ N_{n_0,n_1,n_2}$ be defined by (\ref{eq:P}). Then
\begin{equation}\label{eq:N0}
 N_{0,n_1,n_2}=0\hspace{5mm}{\rm for\;all\;}\; n_1,n_2\geq 0
\end{equation}
\end{Proposition}

\noindent
{\it Proof. } Note that modulo $u_0$
\begin{eqnarray*}
\lefteqn{u_0^{-J_0} u_1^{-J_1} u_2^{-J_2}\Omega(u,J)\equiv} \\
&& \equiv
\left(\sum_{m_1\geq 0}
\frac{(1+ J_0+3 J_1)_{3m_1}}{
{(1+ J_1)_{m_1}}\!^3} u_1^{m_1}\right)
\left(\sum_{m_2\geq 0}
\frac{(1+ J_0+3 J_2)_{3m_2}}{
{(1+ J_2)_{m_2}}\!^3} u_2^{m_2}\right)
\end{eqnarray*}
and take logarithms. The logarithms involve no mixed terms $J_1J_2$.
This shows
$P\equiv 0\bmod u_0$.
{}\hfill$\Box$

\

\

As explained in \cite{F,Sti} a theorem of Bryant and Griffiths
shows
$$
P_j=-{1\over2}U_j\frac{\partial P}{\partial U_j}
$$
for $j=0,1,2$. Hence
\begin{equation}\label{eq:dP}
P_j= \sum_{n_0,n_1,n_2 \geq 0} n_j N_{n_0,n_1,n_2}
\Li_2(U_0^{n_0}U_1^{n_1}U_2^{n_2})
\end{equation}
where $\Li_2(x):=\sum_{n\geq 1}\frac{x^n}{n^2}$ is the dilogarithm
function.

It follows from (\ref{eq:dP}) and (\ref{eq:N0}) that we can get all
numbers
$N_{n_0,n_1,n_2}$ from $P_0$.
The computations are now greatly simplified by observing:

\begin{Lemma}
In ${\cal  R}_{toric}$ the intersection of the $\Z$-module with
basis
\\ $\{1,J_1,J_2,J_1J_2\}$ and the ideal generated by
$J_0,J_1^2,J_2^2$
is $0\,.$
{}\hfill$\Box$
\end{Lemma}

So for studying $\log U_1\,,$ $\log U_2$ and
$P_0$ we may reduce modulo the ideal $(J_0,J_1^2,J_2^2)\,;$
i.e. replace ${\cal  R}_{toric}$ by
$\Z[\bar J_1,  \bar J_2]/(\bar J_1^2,\bar J_2^2)$.
{From} now on we use $J_1$ resp. $J_2$ to denote the classes of
$\bar J_1$ resp. $\bar J_2$ in the latter ring; so we have in
particular
from now on
$$
J_1^2=J_2^2=0
$$
Let
$$
\tilde\Omega(u,J_1,J_2):=\sum_{m_0,m_1,m_2\geq 0}
\frac{(1+ 3 J_1)_{m_0+3m_1}\cdot
(1+ 3 J_2)_{m_0+3m_2}}{
{m_0}!\,^2\cdot {(1+ J_1)_{m_1}}\!^3\cdot
{(1+ J_2)_{m_2}}\!^3}\:u_0^{m_0} u_1^{m_1} u_2^{m_2}
$$
Then clearly
\begin{equation}\label{eq:logtildomexpand}
\begin{array}{rcl}
\log\tilde\Omega(u,J_1,J_2)&=& \log\Omega^{(0)}(u)\, + \,
(\log U_1-\log u_1)\,J_1\,+\\[.4em]
&& +\,(\log U_2-\log u_2)\,J_2
\,+\,{1\over9}P_0\,J_1J_2
\end{array}
\end{equation}
We have the following expansion of $\tilde\Omega(u,J_1,J_2)$ w.r.t.
$u_0$
$$
\tilde\Omega(u,J_1,J_2)=\phi_0(u_1,J_1)\phi_0(u_2,J_2) +
\phi_1(u_1,J_1)\phi_1(u_2,J_2)u_0 + {\cal O}(u_0^2),
$$
where we define
\begin{eqnarray}
\label{eq:phi0}
\phi_0(w,\rho)&:=&\sum_{n\geq0}
{(1+3\rho)_{3n} \over
 {(1+\rho)_n}^3 } w^n \;, \\
\phi_1(w,\rho)&:=&\sum_{n\geq0}
{(1+3\rho)_{1+3n} \over
  {(1+\rho)_n}^3 } w^n  .
\nonumber
\end{eqnarray}
Note
$$
\phi_1(w,\rho)=(1+3\rho)\phi_0(w,\rho)+
3w{\partial\over\partial w}\phi_0(w,\rho)\,.
$$
This shows that modulo $u_0^2$
\begin{equation}\label{eq:lophi}
\displaystyle{
\log\tilde\Omega(u,J_1,J_2)} \equiv
 \displaystyle{\log\phi_0(u_1,J_1)\,+\,\log\phi_0(u_2,J_2)\,+
\hspace{2cm}}
\end{equation}
$$
+ \displaystyle{
\left(1+3J_1+3u_1{\partial\over\partial
u_1}\log\phi_0(u_1,J_1)\right)
\left(1+3J_2+3u_2{\partial\over\partial
u_2}\log\phi_0(u_2,J_2)\right)
u_0}
$$
Comparing (\ref{eq:logtildomexpand}) and (\ref{eq:lophi}) we see
that we have
proved:

\begin{Proposition} Define for $j=1,2$ the function $\bar U_j$ by
$$
\log \bar U_j:=\log
u_j\,+\,{\partial\over\partial\rho}\log\phi_0(u_j,\rho)
|_{\rho=0}\,.
$$
Then
\begin{eqnarray}
\nonumber
\log U_j&=&\log \bar U_j\,+\, {\cal O}(u_0)
\\
\label{eq:P0}
{1\over9}P_0&=&9\left(u_1{\partial\over\partial u_1}\log \bar
U_1\right)
\left(u_2{\partial\over\partial u_2}\log \bar U_2\right)
u_0\,+\, {\cal O}(u_0^2)
\end{eqnarray}
{}\hfill$\Box$
\end{Proposition}

\

Before we can draw conclusions for the numbers $N_{1,n_1,n_2}$ we
must first
analyse $U_0$ modulo $u_0^2$. Let
$$
\tilde{\tilde\Omega}(u,J_0):=\sum_{m_0,m_1,m_2\geq 0}
\frac{(1+ J_0)_{m_0+3m_1}\cdot
(1+ J_0)_{m_0+3m_2}}{
{(1+ J_0)_{m_0}}\!^2\cdot m_1!\,^3\cdot
m_2!\,^3}\:u_0^{m_0} u_1^{m_1} u_2^{m_2}
$$
with as before $J_0^2=0$. Then
$$
\log\tilde{\tilde\Omega}(u,J_0)\,=\,\log\Omega^{(0)}(u)+
(\log U_0-\log u_0)\,J_0
$$
Let
\begin{equation}\label{eq:xi}
\xi(w,\rho):=\sum_{n\geq0}
{(1+\rho)_{3n} \over n!\,^3 } w^n \;
\end{equation}
Then
$$
\tilde{\tilde\Omega}(u,J_0)= \xi (u_1,J_0)\cdot\xi (u_2,J_0)
\,+\, {\cal O}(u_0)
$$
and hence
\begin{equation}\label{eq:U0}
U_0=u_0\cdot\psi(u_1)\cdot\psi(u_2)\,+\, {\cal O}(u_0^2)
\end{equation}
where
$$
\psi(w):=\exp({\partial\over\partial\rho}\log\xi(w,\rho)|_{\rho=0})
$$
By combining (\ref{eq:dP}), (\ref{eq:P0}) and (\ref{eq:U0}) we find

\begin{Corollary}
$$
\sum_{n_1,n_2 \geq 0} N_{1,n_1,n_2}
\bar U_1^{n_1}\bar U_2^{n_2}\,=\,81
\left({1\over\psi(u_1)}u_1{\partial\over\partial u_1}\log \bar
U_1\right)
\left({1\over\psi(u_2)}u_2{\partial\over\partial u_2}\log \bar
U_2\right)
$$
The number $N_{1,n_1,n_2}$ factorizes as
$$
N_{1,n_1,n_2}=b_{n_1} b_{n_2} \;\;,
$$
where the numbers $b_n$ are defined by
\begin{eqnarray}
\sum_{n\geq 0} b_n \bar U_1^n \::=\:9
\left({1\over\psi(u_1)}u_1{\partial\over\partial u_1}\log \bar
U_1\right).
\label{eq:b-prep}
\end{eqnarray}
{}\hfill$\Box$
\end{Corollary}

\begin{Corollary}  Let $\{ b_n \}$ be the sequence of integers
defined by the
expansion $(\ref{eq:b-prep})$.
We obtain the asymptotic expansion of the B-model prepotential as
follows:
 \begin{eqnarray}
 \Psi_B(U_0, U_1, U_2) = \mbox{topological term} +  U_0 B(t_1)
B(t_2) +  {\cal
 O}(U_0^2)
 \label{eq:b-asymp}
 \end{eqnarray}
 where $B(t)$ is defined by the series
 $$
 B(t) = \sum_{n \geq 0} b_n \exp(2 \pi i n t) = \sum_{n \geq 0} b_n
U^n.
 $$
\end{Corollary}

{From} the asymptotic expansions of (\ref{eq:2-A-prep}) and
(\ref{eq:b-asymp}) we
obtain the
following precise identity between two functions, which  actually
follows
from the  Mirror Conjecture~\ref{c:mirror}.

\begin{Conjecture}\label{c:mirror-2}
We will obtain the following identity
$$
\fbox{$ A^{res}(t) \equiv B(t)$}
$$
or equivalently
$$
\fbox{$\sum_{n\geq 0} b_n U^n =
\Theta_{E_8}(3t,t\gamma)\prod_{n\geq 1}(1- U^{3n})^{-12}$}
$$
where $U = \exp(2 \pi i t) $ and $\gamma=(1,1,1,1,1,1,1,-1)$.
\end{Conjecture}

Unfortunately we are unable to prove Conjecture~\ref{c:mirror-2}.
However
since we can explicitly expand the right hand side of
(\ref{eq:b-prep}), we can
obtain
the expansion of $B(t)$ by using a computer and compare the result
with {the expansion} of $A^{res}(t)$.

\begin{Proposition}
Conjecture~\ref{c:mirror-2}  is true up to  order  $U^{50}$.
\end{Proposition}

\vskp
To get started on the computer one may notice:
\begin{eqnarray}
\phi_0(u,0)=\xi(u,0)&=&\sum_{n\geq0}
{(3n)!\over n!\,^3 } u^n \\
{\partial\over\partial u}\phi_0(u,\rho)|_{\rho=0}&=&
\sum_{n\geq0}
{(3n)!\over n!\,^3 } 3(g(3n)-g(n))u^n\\
{\partial\over\partial u}\xi(u,\rho)|_{\rho=0}&=&
\sum_{n\geq0}
{(3n)!\over n!\,^3 } g(3n) u^n
\end{eqnarray}
where
$$
g(n)=\sum_{k=1}^n{1\over k}\;,\hspace{10mm}
g(3n)=\sum_{k=1}^{3n}{1\over k}
$$
A simple PARI program then yields:
\begin{eqnarray}
&& \quad  B(t) =
  9{1\over\psi(u)}u{\partial\over\partial u}\log U =  \nonumber \\
&&9+ 36U+ 126{U^2}+ 360{U^3}+ 945{U^4}+
2268{U^5}+ 5166{U^6}+ 11160{U^7} \nonumber \\
&& + 23220{U^8} +46620{U^9}+ 90972{U^{10}} +
172872{U^{11}} + 321237{U^{12}} \nonumber \\
&&+ 584640{U^{13}}+ 1044810{U^{14}}+1835856{U^{15}}+
3177153{U^{16}}+5421132{U^{17}}
\nonumber \\
&& + 9131220{U^{18}} +15195600{U^{19}}+25006653{U^{20}} +
40722840{U^{21}} \nonumber \\
&& + 65670768{U^{22}}+104930280{U^{23}}+166214205{U^{24}}+
261141300{U^{25}}  \nonumber \\
&&+ 407118726{U^{26}}+  630048384{U^{27}}  + 968272605{U^{28}} +
1478208420{U^{29}}
\nonumber \\
 &&  + 2242463580{U^{30}}+ 3381344280{U^{31}}  +
5069259342{U^{32}} + 7557818940{U^{33}} \nonumber \\
&& + 11208455370{U^{34}}+16538048640{U^{35}}   +
24282822798{U^{36}} \nonumber \\
&& +35487134928{U^{37}}
 + 51626878470{U^{38}}+   74779896240{U^{39}} \nonumber \\
&& +   107861179482{U^{40}}  + 154945739844{U^{41}}   +
221711362038{U^{42}} \nonumber \\
&&+316042958880{U^{43}}  +448856366490{U^{44}}+ 635216766732{U^{45}}
 \nonumber \\
&&+895854679650{U^{46}} + 1259213600736{U^{47}} +
1764210946995{U^{48}}  \nonumber \\
&&  +2463949037340{U^{49}}+3430694064888{U^{50}}  + O(U^{51})
\nonumber \\
&&  \label{eq:b-expand}
\end{eqnarray}

Comparing  this expansion~(\ref{eq:b-expand}) with Table 2 in
Section~\ref{s:rest},
we see that $a_n=  b_n$ for $n\leq 50$.
{}\hfill$\Box$

\section{Appendix I: B-model equation}
\label{s:B-model equations}
In this appendix we derive the equations stated in proposition
\ref{orbifold}
for the mirror $W^*$ of Schoen's example $W$. We use the mirror
construction of
Batyrev-Borisov \cite{Batyrev-Borisov} by means of reflexive
Gorenstein cones
of index $2$. As explained in \cite{Sti} the story in
\cite{Batyrev-Borisov} about split Gorenstein cones and NEF
partitions
can for examples like $W$ be reformulated in terms of triangulations
of the polytope $\Delta$ on the mirror side; more specifically,
$W$ can be embedded in
$\P1\times \P2\times \P2$ as the intersection of a hypersurface of
degree $(1,3,0)$ and a hypersurface of degree $(1,0,3)$.
This leads to the matrix $\sB$ in (\ref{eq:Bmat}) and to the set
$\{3,4\}\times\{5,6,7\}\times\{8,9,10\}$ in (\ref{eq:union}).
To get the reflexive Gorenstein cone $\Lambda$ from which the mirror
of
Schoen's example can be constructed one should take a $7\times 10$
-matrix
$\sA=(a_{ij})$ with rank $7$ and with integer entries such that
$\sA\cdot\sB^t\:=\:0\,.$ We take
$$
\sA\::=\:\left(
\begin{array}{rrrrrrrrrr}
1&0&1&0&1&1&1&0&0&0\\
0&1&0&1&0&0&0&1&1&1\\
0&0&1&-1&0&0&0&0&0&0\\
0&0&0&0&1&-1&0&0&0&0\\
0&0&0&0&1&0&-1&0&0&0\\
0&0&0&0&0&0&0&1&-1&0\\
0&0&0&0&0&0&0&1&0&-1
\end{array}
\right)
$$
Let $\ga_1,\ldots,\ga_{10}\in\Z^7$ be the columns of $\sA$.
Then
$$
\Lambda:=\R_{\geq 0}\ga_1+\ldots+\R_{\geq 0}\ga_{10}\;\subset\R^7
$$
The polytope
$\Delta$ is the convex hull of the points
$\ga_1,\ldots,\ga_{10}$ in $\R^7$. With
$\{3,4\}\times\{5,6,7\}\times\{8,9,10\}$ viewed
as a collection of subsets of $\{1,\ldots,10\}$ the complements of
these
18 subsets are the index sets for the maximal simplices in a
triangulation of
$\Delta$.

Let $\cS_\Lambda$ denote the subalgebra of the algebra of Laurent
polynomials
$$
\C[u_1^{\pm 1},\ldots,u_7^{\pm 1}]
$$ generated by the monomial
$u_1^{m_1}\cdot\ldots\cdot u_7^{m_7}$ with
$(m_1,\ldots,m_7)\in\Lambda\cap\Z^7$. Giving such a monomial degree
$m_1+m_2$
makes $\cS_\Lambda$ a graded ring.
The scheme
$\BP_\Lambda\::=\:{\rm Proj }\cS_\Lambda$ is a projective toric
variety
of dimension $6$.
A global section of $\cO_{\BP_\Lambda}(1)$ is given by a Laurent
polynomial (with coefficients $v_1,\ldots,v_{10}$)
\begin{eqnarray*}
\es & = & u_1(v_1+v_5u_4u_5+v_6u_4^{-1}+v_7u_5^{-1}+v_3u_3) + \\
& & u_2(v_2+v_8u_6u_7+v_9u_6^{-1}+v_{10}u_7^{-1}+v_4u_3^{-1})
\end{eqnarray*}
For generic coefficients $v_1,\ldots,v_{10}$ the
zero locus of $\es$ in $\BP_\Lambda$ is a generalized Calabi-Yau
manifold of
dimension $5$ in the sense of
\cite{Batyrev-Borisov}. This is one mirror of $W$ suggested by
\cite{Batyrev-Borisov}.

As in \cite{Batyrev-Borisov} Section 4,  one can also realize a
mirror as a
complete intersection Calabi-Yau threefold in a $5$-dimensional
toric variety,
as follows.
$\BP_\Lambda$ is a compactification of the torus
$(\C^*)^7/\C^*$ where $\C^*:=\{(u,u,1,1,1,1,1)\in (\C^*)^7\}$.
The morphism
\begin{eqnarray*}
(\C^*)^7&\rightarrow&
\P1\times\P1\times\P3\times\P3
\end{eqnarray*}
given by
\begin{eqnarray*}
\lefteqn{(u_1,\ldots,u_7)\mapsto ([u_1:u_1u_3],[u_2:u_2u_3^{-1}], }
\\
&& [u_1:u_1u_4u_5:u_1u_4^{-1}:u_1u_5^{-1}],
[u_2:u_2u_6u_7:u_2u_6^{-1}:u_2u_7^{-1}])
\end{eqnarray*}
extends to a morphism
$\BP_\Lambda\rightarrow \P1\times\P1\times\P3\times\P3$. The image
is
$$
V:=\left\{\left.
\begin{array}{l}
[p_0:p_1]\times [q_0:q_1]\times [s_0:s_1:s_2:s_3]\times
[t_0:t_1:t_2:t_3]  \\
 \in \P1\times\P1\times\P3\times\P3 \\
p_0q_0= p_1q_1,  \quad s_0^3= s_1s_2s_3, \quad
t_0^3 = t_1t_2t_3
\end{array}
\right.\right\}
$$
As noted in \cite{Batyrev-Borisov} Cor.3.4 the complement of the
generalized
Calabi-Yau $5$-fold $\es=0$
in $\BP_\Lambda$ is a
$\C$-bundle over the complement in $V$ of the complete intersection
Calabi-Yau $3$-fold with equations
\begin{eqnarray*}
(v_1s_0+v_5s_1+v_6s_2+v_7s_3)p_0+v_3s_0p_1 &=&0\\
(v_2t_0+v_8t_1+v_9t_2+v_{10}t_3)q_0+v_4t_0q_1&=&0
\end{eqnarray*}
This complete intersection Calabi-Yau $3$-fold itself is another
realization for a mirror of $W$.

Now note that the morphism
\begin{eqnarray*}
\P1\times\P2\times\P2&\rightarrow&\P1\times\P1\times\P3\times\P3\;,
\end{eqnarray*}
$$
\begin{array}{l}
([z_0:z_1],[x_0:x_1:x_2],[y_0:y_1:y_2]) \mapsto \\
([z_0:z_1],[z_1:z_0],[x_0x_1x_2:x_0^3:x_1^3:x_2^3]
,[y_0y_1y_2:y_0^3:y_1^3:y_2^3])
\end{array}
$$
realizes $V$ also as the quotient of $ \P1 \times \P2 \times \P2$
by the group $ \Z_3\times\Z_3 $ acting as in
Proposition~\ref{orbifold}.
This completes the proof of Proposition~\ref{orbifold}.

\section{Appendix II: The Theta function of the $E_8$ lattice }
 \label{s:appendixII}

Let $\Lambda $ be a lattice of rank $d$ with positive definite
quadratic form $Q:\Lambda \lra \Z$.
We can fix an embedding  $\Lambda \hookrightarrow \R^d$
such that the quadratic form $Q$
is induced by the usual Euclidean inner product $( \: , \: )$.
Let $\cH = \{ \tau \in \C | \mbox{\rm Im} (\tau) > 0 \}$ be the
upper
half plane. We denote by $\w = (w_1, \cdots, w_d)$ the
standard complex coordinates of $\C^d =  \R^d \otimes \C$.
We define the  theta function associated to the lattice $\Lambda$ by
\begin{eqnarray}
\Theta_{\Lambda}(\tau, \w) =
\sum_{\sigma \in \Lambda}
\exp(2 \pi i ((\tau/2) Q(\sigma) +   (\sigma, \w)).
\end{eqnarray}
For  certain Calabi-Yau 3-folds with a fibration by
abelian surfaces
one can calculate a part of the prepotential of the Yukawa coupling
arising from the
sections of the fibration by using the theta function associated to
the Mordell-Weil lattice  \cite{Saito}. Since the Mordell-Weil
lattice of
a generic Schoen's example is isometric to $E_8 \times E_8$, we
would like to
calculate the theta function of  $E_8$ and write it in an explicit
form.

For that purpose we fix a standard embedding of $D_8$ and $E_8$
into $\R^8$. (\cite{C-S} p. 117 $\sim$
p. 121).
Let $e_1, e_2, \cdots, e_8$ be the
standard orthonormal basis of $\R^8$. An element of $\R^8$
is written as $\sum_{i=1}^{8} x_i e_i$. We define  lattices in
$\R^8$
$$
\Z^{8} := \left\{ \sum_{i=1}^{8} x_i e_i,  x_i \in \Z \right\}
\supset
 D_8 := \left\{  \sum_{i=1}^{8} x_i e_i \in \Z^8, \\
\sum_{i=1}^{8} x_i \equiv 0\, (2) \right\},
$$

$$
E_8 = D_8 \cup (D_8 + s_0),  \quad s_0 = \frac{1}{2}
\sum_{i=1}^{8} e_i,
$$
The inner product $( \:, \:)$ induces positive
definite  bilinear forms on these lattices and $E_8$ and $D_8$ have
integral bases whose intersection matrices are the Cartan matrices
of
$E_8$ and $D_8$ respectively.

The theta function for the one dimensional lattice $\Lambda = \Z$
with $Q(n) = n^2$ is  the Jacobi theta function:
\begin{eqnarray}
\vartheta(\tau, w)  :=\Theta_{\Z}(\tau, w) =   \sum_{n \in \Z}
\exp(\pi i n^2\tau + 2 \pi i n w).
\end{eqnarray}
We also have  the following 4 theta functions (cf.
\cite{Mum-Tata-1}):

\begin{eqnarray}
\vartheta_{0,0}(\tau, w) & = &  \vartheta(\tau, w) \\
\vartheta_{0,1}(\tau, w) & = & \vartheta(\tau, w+\frac{1}{2}) \\
\vartheta_{1,0}(\tau, w) & = & \exp(\frac{\pi i \tau}{4} + \pi i w)
\cdot \vartheta(\tau, w+\frac{\tau}{2}) \\
\vartheta_{1,1}(\tau, w) & = & \exp(\frac{\pi i \tau}{4} + \pi
i(w+\frac{1}{2}))
\cdot \vartheta(\tau, w+\frac{\tau+1}{2})
\end{eqnarray}

\begin{Proposition} \label{pr:e8theta}
 Let $\w = (w_1, w_2, \cdots, w_8) \in \C^8$.
\begin{eqnarray}
\Theta_{\Z^8}(\tau, \w) & =  & \prod_{i=1}^8
\vartheta_{0,0}(\tau, w_i)  \label{eq:z8}\\
\Theta_{E_8}(\tau, \w) &= & \frac{1}{2} \sum_{(a, b) \in
(\Z/2\Z)^2}
\prod_{i=1}^8 \vartheta_{a,b}(\tau, w_i) \label{eq:e8}
\end{eqnarray}
\end{Proposition}

\noindent
{\it Proof.} Straightforward exercise. See also \cite{D-G-W}.
{}\hfill$\Box$

\

Recall $\gamma=(1,1,1,1,1,1,1,-1)$. The above formulas show
(cf.\cite{Mum-Tata-1}):
\begin{eqnarray*}
\vartheta_{0,0}(\tau,-w)=\vartheta_{0,0}(\tau,w)\:,&\hspace{3mm}&
\vartheta_{0,1}(\tau,-w)=\vartheta_{0,1}(\tau,w)\:,\\
\vartheta_{1,0}(\tau,-w)=\vartheta_{1,0}(\tau,w)\:,&\hspace{3mm}&
\vartheta_{1,1}(\tau,-w)=-\vartheta_{1,1}(\tau,w)
\end{eqnarray*}
and hence
\begin{equation}
\Theta_{E_8}(3t, t\gamma)=\frac{1}{2}\{ \vartheta_{0,0}(3t,t)^8
+\vartheta_{0,1}(3t,t)^8+\vartheta_{1,0}(3t,t)^8
-\vartheta_{1,1}(3t,t)^8 \}
\end{equation}
Next note:
\begin{eqnarray*}
\vartheta_{0,0}(3t,t)&=&
\exp(-\pi i t/3)\sum_{ m\equiv\pm 1\,(3)}\exp(\pi i t m^2/3) \\
\vartheta_{0,1}(3t,t)&=&
-\exp(-\pi i t/3)\sum_{ m\equiv\pm 1\,(3)}(-1)^m\exp(\pi i t m^2/3)
\\
\vartheta_{1,0}(3t,t)&=&
\exp(-\pi i t/3)
\sum_{ m\equiv\pm 1\,(6)}\exp(\pi i t m^2/12) \\
\vartheta_{1,1}(3t,t)&=&
-i\exp(-\pi i t/3)
\sum_{ m\equiv\pm 1\,(6)} \chi (m)\,\exp(\pi i t m^2/12)
\end{eqnarray*}
where the summations run over $m \in \N$ with the indicated
restrictions
and $\chi (m)=1$ (resp. $=-1$) if $m\equiv\pm 1\bmod 12$
(resp. $\equiv\pm 5\bmod 12$). Another useful observation is that
the Jacobi product formula for $\vartheta_{1,1}(\tau,w)$ (see
\cite{Mum-Tata-1})
implies
$$
\vartheta_{1,1}(3t,t)=-i\exp(-\pi i t/4)\prod_{m\geq 1}(1-\exp(2\pi
i m t))
$$

Now the computer can do its work and compute the expansion of
$\Theta_{E_8}(3t, t\gamma)$.

 \section*{Acknowledgments}

The first author  would like to thank J.Bryan and N.C.Leung for
notifying him of the paper\cite{G-P}.
He would like to thank also to S.-T.Yau and the
Mathematics Department
of Harvard University for their hospitality when finishing this
work.

The second author  would like to thank
Taniguchi foundation for their generous support for the Symposium.
He would like to thank
also all participants in the symposium  with whom he enjoyed
fruitful
discussion.  In particular, He would like to thank Ron Donagi for the
discussion about \cite{D-G-W}.
He would like to thank the staff of Kobe University, where he is
enjoying daily
stimulating atmosphere and discussion about mathematics. Special thanks
are due to
 Kota Yoshioka in Kobe University who kindly remarked
Lemma~\ref{l:hilb-2}.

The third author would like to thank the Japan Society for
the Promotion of Science for a JSPS Invitation Fellowship
in November-December 1996 and Kobe University
for support for a visit in July 1997.
He expresses special thanks to his host, Masa-Hiko Saito, for
creating
a very stimulating atmosphere during these two visits to Kobe.

\vskp
\vskp


\end{document}